\documentclass[11pt]{article}
\usepackage{amsfonts}
\usepackage{latexsym}
\usepackage{amssymb}
\newcommand{\cF}{{\cal F}}

\newcommand{\cI}{{\cal I}}
\newcommand{\cL}{{\cal L}}

\newcommand{\cP}{{\cal P}}

\newcommand{\cS}{{\cal S}}

\newcommand{\mcH}{\mathcal{H}}

\newcommand{\cV}{{\cal V}}

\newcommand{\ch}{\mathfrak{h}}

\renewcommand{\H}{\mathbb{H}}

\newcommand{\R}{\mathbb{R}}

\renewcommand{\H}{\mathbb{H}}

\newcommand{\gm}{\gamma}

\newcommand{\lb}{\lambda}
\newcommand{\Lb}{\Lambda}
\newcommand{\ep}{\varepsilon}

\newcommand{\ph}{\varphi}

\newcommand{\sm}{\setminus}

\newcommand{\lan}{\langle}
\newcommand{\ran}{\rangle}
\newcommand{\lls}{\mbox{\large $($}}
\newcommand{\rls}{\mbox{\large $)$}}
\newcommand{\lLs}{\mbox{\Large $($}}
\newcommand{\rLs}{\mbox{\Large $)$}}

\newcommand{\ra}{\rightarrow}
\newcommand{\lra}{\longrightarrow}

\newcommand{\cir}{{\scriptstyle \circ}}
\newcommand{\der}{\partial}

\newcommand{\avint}{\hbox{\vrule height3.5pt depth-2.8pt
width4pt}\mkern-12mu\int\nolimits}
\newcommand{\diam}{\mbox{diam}}

\newcommand{\Lip}{\mbox{Lip}}

\renewcommand{\span}{\mbox{span}}

\newcommand{\im}{\mbox{Im$\,$}}
\renewcommand{\ker}{\mbox{Ker$\,$}}
\newcommand{\vol}{\mbox{{\rm vol}}}

\newcommand{\Btil}{\tilde{B}}

\newcommand{\Ftil}{\tilde{F}}

\newcommand{\gtil}{\tilde{g}}

\newcommand{\xtil}{\tilde{x}}
\newcommand{\ytil}{\tilde{y}}

\newtheorem{The}{Theorem}[section]
\newtheorem{Lem}[The]{Lemma}
\newtheorem{Def}[The]{Definition}
\newtheorem{Rem}[The]{Remark}
\newtheorem{Pro}[The]{Proposition}
\newtheorem{Cor}[The]{Corollary}
\newtheorem{Exa}[The]{Example}

\begin{document}
\title{{\bf Blow-up of regular submanifolds in Heisenberg groups and
applications}}
\author{{\sc Valentino Magnani}
\\ Dipartimento di Matematica\\ via Buonarroti n.2, 56127, Pisa, Italy \\
e-mail magnani@dm.unipi.it}
\date{}
\maketitle

\begin{abstract}
We obtain a blow-up theorem for regular submanifolds in the Heisenberg group, where intrinsic dilations are used. Main consequence of this result is an explicit formula for the density of $(p$$+$$1)$-dimensional spherical Hausdorff measure restricted to a $p$-dimensional submanifold with respect to the Riemannian surface measure.
We explicitly compute this formula in some simple examples and we present a lower semicontinuity result for the spherical Hausdorff
measure with respect to the weak convergence of currents.
Another application is the proof of an intrinsic coarea formula for vector-valued mappings on the Heisenberg group.  
\end{abstract}

\vskip1.5truecm
\tableofcontents

\pagebreak

\section{Introduction}

In recent years, several efforts have been devoted to the
project of developing Analysis and Geometry in stratified groups and
more general Carnot-Carath\'eodory spaces with several monographs
and surveys on this subject. Among them we mention
\cite{BelRis}, \cite{FS}, \cite{HajKos}, \cite{Kup}, \cite{Montgom}, \cite{VSC}, but this list could be surely enlarged.

Our study fits into the recent project of developing Geometric
Measure Theory in these spaces.
Ambient of our investigations is the $(2n$$+$$1)$-dimensional Heisenberg group $\H^n$, which represents the simplest model
of non-Abelian stratified group, \cite{FS}, \cite{Stein}. 
Aim of this paper is to present an intrinsic blow-up theorem for
$C^1$ sub\-manifolds in the geo\-metry of the Heisenberg group
along with its applications.
The main feature of this procedure is the use of natural
dilations of the group, namely, a one-parameter family of 
group homomorphisms that are homogeneous with respect to
the distance of the group.
Recall that dilations in $\H^n$ are anisotropic,
hence they differently act on different directions of
the submanifold. The foremost directions are the so-called
horizontal directions, that determine the ``sub-Riemannian geometry" of the Heisenberg group: at any point $x\in\H^n$
a $2$$n$-dimensional subspace $H_x\H^n\subset T_x\H^{2n+1}$ is given
and the family of all horizontal spaces $H_x\H^n$ forms the
so-called horizontal subbundle $H\H^n$.
We will defer full definitions to Section~\ref{basic}.

The blow-up procedure consists in enlarging the submanifold
$\Sigma$ at some point $x\in\Sigma$ by intrinsic dilations and taking the intersection of the magnified submanifold with a bounded set
centered at $x$.
We are interested in studying the case when $T_x\Sigma\not\subset H_x\H^n$, namely, $x$ is a transverse point. The effect of rescaling the submanifold at a transverse point $x$ can be
obtained by considering the behavior of
$\vol_p(B_{x,r}\cap\Sigma)/r^{p+1}$ as $r\ra0^+$, that heuristically is
\begin{eqnarray*}
\frac{\vol_p(B_{x,r}\cap\Sigma)}{r^{p+1}}=\frac{
\vol_p\lls l_x\delta_r(B_1\cap\Sigma_{x,r})\rls}{r^{p+1}}=\frac{
\vol_p\lls \delta_r(B_1\cap\Sigma_{x,r})\rls}{r^{p+1}}\approx
\alpha(x)\,\vol_p(B_1\cap\Sigma_{x,r}).
\end{eqnarray*}
Here $\vol_p$ denotes the $p$-dimensional Riemannian
measure restricted to $\Sigma$,
the left translation $l_x:\H^n\lra\H^n$ is given by $l_x(y)=x\cdot y$,
the dilation of factor $r>0$ is $\delta_r:\H^n\lra\H^n$,
the dilated submanifold at $x$ is $\Sigma_{x,r}=\delta_{1/r}\lls l_{x^{-1}}\Sigma\rls$ and $B_{x,r}$ is the open ball of center $x$ and radius $r$ with respect to a fixed homogeneous distance.
The meaning of $\alpha(x)$ will be clear in the following theorem,
that makes rigorous our previous consideration and represents
our first main result.
%
%
%
%%%%%%%%%%%%%%%%%%%%%%%%%%%%%%%%%%%%%%%%%%%%%%%%%%%%%%%%%%%%%
%         BLOW-UP OF CODIMENSION k SUBMANIFOLDS
%%%%%%%%%%%%%%%%%%%%%%%%%%%%%%%%%%%%%%%%%%%%%%%%%%%%%%%%%%%%%
%
%
%
\begin{The}[Blow-up]\label{bwpt}
Let $\Sigma$ be a $p$-dimensional $C^1$ submanifold of $\Omega$,
where $\Omega$ is an open subset of $\H^n$ and let $x$
be a transverse point. Then the following limit holds
\begin{eqnarray}\label{keyblow}
\lim_{r\ra0^+}\frac{\vol_p(\Sigma\cap B_{x,r})}{r^{p+1}}
=\frac{\theta_p^\rho\lls\tau_{\Sigma,\cV}(x)\rls}
{|\tau_{\Sigma,\cV}(x)|}.
\end{eqnarray}
\end{The}
A novel object appearing in this limit is the vertical tangent $p$-vector $\tau_{\Sigma,\cV}(x)$, introduced in Definition~\ref{vertangpvect}.
Its associated $p$-dimensional subspace of $\ch^{n}$ is
a subalgebra whose image through the exponential map
represents the blow-up limit of the rescaled submanifold
$\Sigma_{x,r}$ as $r\ra0^+$.
The $p$-vector $\tau_{\Sigma,\cV}(x)$ in higher codimension
plays the same role that the well known horizontal normal
$\nu_H$ plays in codimension one (compare for instance with \cite{Mag2}). 
The metric factor $\theta\lls\tau_{\Sigma,\cV}(x)\rls$,
introduced in \cite{Mag2}, corresponds to the measure of the intersection of $B_1$ with the vertical subspace associated to the vertical tangent $p$-vector $\tau_{\Sigma,\cV}(x)$.
A first consequence of Theorem~\ref{bwpt} is an explicit
formula to compute the $(p$$+$$1)$-dimen\-sional spherical Hausdorff measure of $p$-dimensional $C^1$ submanifolds in the Heisenberg group.
In fact, thanks to $\cS^{p+1}$-negligibility of characteristic
points proved in \cite{Mag5}, Theorem~\ref{bwpt} along with standard theorems on differentation of measures, immediately give
the following result.
%
%
%
%%%%%%%%%%%%%%%%%%%%%%%%%%%%%%%%%%%%%%%%%%%%%%%%%%%%%%%%%%%%%%
% SPHERICAL HAUSDORFF MEASURE OF SUBMANIFOLDS CONSTRUCTED BY
% HOMOGENEOUS DISTANCE WITH CONSTANT METRIC FACTOR
%%%%%%%%%%%%%%%%%%%%%%%%%%%%%%%%%%%%%%%%%%%%%%%%%%%%%%%%%%%%%%
%
%
%
\begin{The}\label{constmetintsp}
Let $\rho$ be a homogeneous distance with constant metric factor $\alpha>0$
and let $\cS^{p+1}_{\H^{n}}=\alpha\,\cS^{p+1}_\rho$. 
Then we have
\begin{eqnarray}\label{thef2c}
\cS^{p+1}_{\H^{n}}(\Sigma)=\int_\Sigma |\tau_{\Sigma,\cV}(x)|\,d\vol_p(x).
\end{eqnarray}
\end{The}
Note that in codimension one, the integral formula (\ref{thef2c}) fits into the results of \cite{Mag2} in stratified groups.
The connection between these results is shown in Proposition~\ref{cod1h}.
There are several examples of homogeneous distances
satisfying hypothesis of Theorem~\ref{constmetintsp}, as we show in
Example~\ref{examphomd}. Proposition~\ref{homconstmet}
shows a class of homogeneous distances having constant metric factor.
Proposition~\ref{concreteex} shows how the computation of the $(p$$+1)$-dimensional spherical Hausdorff measure of a submanifold
can be easily performed in several examples, that will appear in Section~\ref{sphersub}.

Another consequence of Theorem~\ref{bwpt} is the validity of
an intrinsic coarea formula for vector-valued Lipschitz mappings
defined on the Heisenberg group. 
By Sard theorem and the classical Whitney approximation theorem we can assume that a.e. level set is a submanifold of class $C^1$, then we apply representation formula (\ref{thef2c}). The core of the proof stands in the key relation
\begin{eqnarray}\label{forint}
|\tau_{\Sigma,\cV}(x)|=\frac{J_Hf(x)}{J_gf(x)}\,,
\end{eqnarray}
which surprisingly connects vertical tangent $p$-vector with horizontal jacobian $J_Hf$.
The proof of (\ref{forint}) is given in Theorem~\ref{jhjg}.
Thus, we can establish the following result.
%
%
%
%%%%%%%%%%%%%%%%%%%%%%%%%%%%%%%%%%%%%%%%%%%%%%%%%%%%%%%%%%%%%%
%                     COAREA FORMULA
%%%%%%%%%%%%%%%%%%%%%%%%%%%%%%%%%%%%%%%%%%%%%%%%%%%%%%%%%%%%%%
%
%
%
\begin{The}[Coarea formula]\label{coacodk}
Let $f:A\lra\R^k$ be a Riemannian Lipschitz map, where
$A\subset\H^{n}$ is a measurable subset and $1\leq k<2n+1$.
Let $\rho$ be a homogeneous distance with constant metric
factor $\alpha>0$. Then for every measurable function
$u:A\lra[0,+\infty]$ the formula
\begin{eqnarray}\label{cht}
\int_A u(x)\,J_{H}f(x)\,dx=\int_{\R^k}\left(\int_{f^{-1}(t)\cap A}u(y)
\;d\cS_{\H^{n}}^{p+1}(y)\right)\,dt
\end{eqnarray}
holds, where $p=2n+1-k$ and $\cS_{\H^n}^{p+1}=\alpha\,\cS_\rho^{p+1}$.
\end{The}
This coarea formula along with that of \cite{Mag4}, which
is a particular case, represent first examples of intrinsic
coarea formulae for vector valued mappings defined on
non-Abelian Carnot groups.
It remains an interesting open question the extension of
coarea formula to Lipschitz mappings with respect to
a homogeneous distance. Only in the case of real-valued mappings
this problem has been settled in \cite{Mag3}.
This question is intimately related to a
blow-up theorem of ``intrinsicly regular" submanifolds.
In this connection, we mention a recent work by Franchi,
Serapioni and Serra Cassano \cite{FSSC6}, where a notion of intrinisic submanifold in $\H^n$ has been introduced in arbitrary codimension.
According the their terminology, a $k$-codimensional $\H$-regular submanifold for algebraic reasons must satisfy
$1\leq k\leq n$. With this restriction it might be highly irregular,
even unrectifiable in the Euclidean sense, \cite{KirSer}.
Nevertheless they show that an area-type formula for its $(p$$+$$1)$-dimensional spherical Hausdorff measure still holds. 
Here we wish to emphasize the difference in our approach, where we consider $C^1$ submanifolds, but with no restriction on their codimension.

Let us summarize the contents of the present paper.
Section~\ref{basic} recalls some notions.
Section~\ref{blwtran} is devoted to the proof of Theorem~\ref{bwpt}.
In Section~\ref{sphersub} we show the validity of Theorem~\ref{constmetintsp}, along with its applications.
Precisely, in Theorem~\ref{invarresc}
we show how a suitable rescaling of the spherical Hausdorff measure yields an intrinsic surface measure only depending on the sub-Riemannian metric,
namely, the restriction of the Riemannian metric to the horizontal subbundle. 
In Proposition~\ref{homconstmet}, we single out a privileged class of 
homogeneous distances having constant metric factor.
We present several explicit computations of $(p$$+$$1)$-dimensional
spherical Hausdorff measure in concrete examples.
As another application of Theorem~\ref{constmetintsp}, we show a lower semicontintuity result for the spherical Hausdorff measure with respect to weak convergence of regular currents. Section~\ref{coareasect} establishes
an intrinsic coarea formula for vector-valued Riemannian Lipschitz
mappings on the Heisenberg group.
\vskip.25truecm
\noindent
{\bf Acknowledgments.}
I wish to thank Bruno Franchi, Raul Serapioni and Francesco Serra
Cassano for pleasant discussions on intrinsic surface area in Heisenberg groups.

%
%
%
%
%
%%%%%%%%%%%%%%%%%%%%%%%%%%%%%%%%%%%%%%%%%%%%%%%%%%%%%%%%%%%%%%
%
%               SECTION ON BASIC NOTIONS
%
%%%%%%%%%%%%%%%%%%%%%%%%%%%%%%%%%%%%%%%%%%%%%%%%%%%%%%%%%%%%%%
%
%
%
%
\section{Some basic notions}\label{basic}

The $(2n$$+$$1)$-dimensional Heisenberg group $\H^{n}$
is a simply connected Lie group whose Lie algebra $\ch^{n}$ is
equipped with a basis $(X_{1},\dots,X_{2n},Z)$
satisfying the bracket relations
\begin{eqnarray}\label{commhei}
[X_{k},X_{k+n}]=2\,Z
\end{eqnarray}
for every $k=1,\ldots,n$.
We will identify the Lie algebra $\ch^{n}$ with the isomorphic Lie
algebra of left invariant vector fields on $\H^{n}$, so that 
any $X_j$ also denotes a left invariant vector field of $\H^{n}$.
In the terminology of Differential Geometry, the basis $(X_{1},\ldots,X_{2n},Z)$ forms a moving frame in $\H^n$.
We will say that $(X_{1},\ldots,X_{2n},Z)$ is our {\em standard frame}.
In particular, $(X_{1},\ldots,X_{2n})$ is a {\em horizontal frame}
and it spans a smooth distribution of $2n$-dimensional hyperplanes, called
{\em horizontal hyperplanes} and denoted by $H_x\H^n$ for every
$x\in\H^n$. The collection of all horizontal hyperplanes forms
the so called {\em horizontal subbundle}, denoted by $H\H^{n}$. 
In the sequel, we will fix the unique left invariant Riemannian metric $g$
such that the standard frame $(X_1,X_2,\ldots,X_{2n},Z)$
forms an orthonormal basis at each point.
\begin{Def}{\rm Every set of left invariant vector fields
$(Y_1,\ldots,Y_{2n})$ spanning the horizontal hyperplane at the unit
element of $\H^n$ will be called {\em horizontal frame}.}
\end{Def}
Recall that the exponential map $\exp:\ch^{n}\lra\H^{n}$ is a diffeomorphism, then it is possible to introduce a system of coordinates in all of $\H^n$. 
%
%%%%%%%%%%%%%%%%%%%%%%%%%%%%%%%%%%%%%%%%%%%%%%%%%%%%%%%%%%%%%
%                    GRADED COORDINATES
%%%%%%%%%%%%%%%%%%%%%%%%%%%%%%%%%%%%%%%%%%%%%%%%%%%%%%%%%%%%%
%
\begin{Def}[Graded coordinates]\label{gradef}
{\rm Let $(Y_1,\ldots,Y_{2n})$ be a horizontal frame
and let $W$ be a non horizontal left invariant vector field.
The frame $(Y_1,\ldots,Y_{2n},W)$ defines a coordinate chart
$F:\R^{2n+1}\lra\H^{n}$ given by
\begin{eqnarray}\label{gco}
F(y)=\exp\bigg(y_{2n+1}W+\sum_{j=1}^{2n}y_{j}Y_{j}\bigg).
\end{eqnarray}
Coordinates defined by (\ref{gco}) are called 
{\em graded coordinates} in the case $W=Z$ and 
{\em standard coordinates} in the case the standard frame $(X_1,\ldots,X_{2n},Z)$ is used. In general we will
say that the coordinates are associated to the frame
$(Y_1,\ldots,Y_{2n},W)$
}\end{Def}
We will assume throughout that a system of standard coordinate
is fixed, if not stated otherwise.
\begin{Rem}{\rm
Note that the horizontal frame $(Y_1,\ldots,Y_{2n})$ of Definition~\ref{gradef}
may not satisfy relations (\ref{commhei}), where $X_i$ are replaced by $Y_i$.
}\end{Rem}
The standard frame with respect to standard coordinates reads as follows
%
%%%%%%%%%%%%%%%%%%%%%%%%%%%%%%%%%%%%%%%%%%%%%%%%%%%%%%%%%%%%%
%
%    VECTOR FIELDS WITH RESPECT TO SYMPLECTIC COORDINATES
%
%%%%%%%%%%%%%%%%%%%%%%%%%%%%%%%%%%%%%%%%%%%%%%%%%%%%%%%%%%%%%
%
\begin{eqnarray}\label{hvf}
\tilde{X}_{k}=\der_{x_{k}}-x_{k+n}\,\der_{x_{2n+1}},\;\,
\tilde{X}_{k+n}=\der_{x_{k+n}}+x_{k}\,\der_{x_{2n+1}}\;\,
\mbox{and}\;\,
\tilde{Z}=\der_{x_{2n+1}}
\end{eqnarray}
and the group operation is given by the following formula
%
%%%%%%%%%%%%%%%%%%%%%%%%%%%%%%%%%%%%%%%%%%%%%%%%%%%%%%%%%%%%%
%
%                    GROUP OPERATION 
%
%%%%%%%%%%%%%%%%%%%%%%%%%%%%%%%%%%%%%%%%%%%%%%%%%%%%%%%%%%%%%
%
\begin{eqnarray}\label{glaw}
x\cdot y
=\bigg(x_{1}+y_{1},\ldots,x_{2n}+y_{2n},x_{2n+1}+y_{2n+1}+
\sum_{j=1}^{n}(x_{k}y_{k+n}-x_{k+n}y_{k})\bigg).
\end{eqnarray}
A natural family of dilations which respects the group 
operation (\ref{glaw}) can be defined as follows
\begin{eqnarray}\label{}
\delta_{r}(x)=(rx_{1},rx_{2},\ldots,rx_{2n},r^{2}x_{2n+1})
\end{eqnarray}
for every $r>0$. In fact, the map $\delta_{r}:\H^{n}\lra\H^{n}$
defined above
is a group homomorphism with respect to the operation (\ref{glaw}).

In contrast with Analysis in Euclidean spaces, where the Euclidean 
distance is the most natural choice, in the Heisenberg group several
distances have been introduced for different purposes.
However, all of them are {\em homogeneous} in the following sense.
If $\rho:\H^n\times\H^n\lra[0,+\infty+[$ is a homogeneous distance, then 
\begin{enumerate}
\item
$\rho$ is a continuous with respect to the topology of $\H^n$,
\item
$\rho(xy,xz)=\rho(y,z)$ for every $x,y,z\in\H^{n}$,
\item
$\rho(\delta_{r}y,\delta_{r}z)=r\,\rho(y,z)$ for every $y,z\in\H^{n}$ 
and every $r>0$.
\end{enumerate}
To simplify notations we write $\rho(x,0)=\rho(x)$, where $0$
denotes either the origin of $\R^{2n+1}$ or the unit element of $\H^{n}$.
The open ball of center $x$ and radius $r>0$ with respect to a homogeneous
distance is denoted by $B_{x,r}$.
The Carnot-Carath\'eodory distance is an important example of homogeneous distance, \cite{Gr1}.
However, all of our computations hold for a general homogeneous distance,
therefore in the sequel $\rho$ will denote a homogeneous distance,
if not stated otherwise.
Note that the Hausdorff dimension of $\H^{n}$ with respect
to any homogeneous distance is $2n+2$.
Next, we recall the notion of Riemannian jacobian.
%
%
%
%%%%%%%%%%%%%%%%%%%%%%%%%%%%%%%%%%%%%%%%%%%%%%%%%%%%%%%%%%%%%
%                   RIEMANNIAN JACOBIAN
%%%%%%%%%%%%%%%%%%%%%%%%%%%%%%%%%%%%%%%%%%%%%%%%%%%%%%%%%%%%%
%
%
%
\begin{Def}[Riemannian jacobian]\label{riemjacobian}{\rm
Let $f:M\lra N$ be a $C^1$ smooth mapping of Riemannian manifolds
and let $x\in M$, where $M$ and $N$ have dimension $d$ and $k$, respectively.
The Riemannian jacobian of $f$ at $x$ is given by
\begin{eqnarray}
J_gf(x)=\|\Lb_k\lls df(x)\rls\|,
\end{eqnarray}
where $\Lb_k\lls df(x)\rls:\Lb_d(T_xM)\lra\Lb_k(T_{f(x)}N)$
is the canonical linear map associated to $df(x):T_xM\lra T_{f(x)}N$.
The norm of $\Lb_k\lls df(x)\rls$ is understood
with respect to the induced scalar products on $\Lb_d(T_xM)$ and $\Lb_k(T_{f(x)}N)$. We recall scalar products of $p$-vectors
in (\ref{pscalarproduct}).
}\end{Def}
To compute the Riemannian jacobian, we fix two orthonormal bases 
$(X_1,\ldots,X_d)$ and $(E_1,\ldots,E_k)$
of $T_xM$ and $T_{f(x)}N$, respectively, and we represent
$df(x)$ with respect to these bases by the matrix
\begin{eqnarray}
\nabla_{X,E}f(x)=\left[\begin{array}{cccccc}
\lan E_1,df(x)(X_1)\ran & \lan E_1,df(x)(X_2)\ran & \ldots & \lan E_1,df(x)(X_d)\ran  \\
\lan E_2,df(x)(X_1)\ran & \lan E_2,df(x)(X_2)\ran & \ldots & \lan E_2,df(x)(X_d)\ran  \\
\vdots        & \vdots        & \vdots & \vdots     \\
\lan E_k,df(x)(X_1)\ran & \lan E_k,df(x)(X_2)\ran & \ldots &  \lan E_k,df(x)(X_d)\ran  \\
\end{array}\right].
\end{eqnarray}
Then the jacobian of the matrix $\nabla_{X,E}f(x)$
coincides with $J_gf(x)$. In the sequel, it will be useful
to fix the following notation to indicate minors of a matrix.
\begin{Def}{\rm
Let $G$ be an $m\times n$ matrix with $m\leq n$. We denote by 
$G_{i_{1}i_{2}\ldots i_{m}}$ the $m\times m$ submatrix with
columns $(i_{1},i_{2},\ldots,i_{m})$. We define the minor
\begin{eqnarray}
M_{i_{1}i_{2}\ldots i_{m}}(G)=\det\left(G_{i_{1}i_{2}\ldots i_{m}}\right).
\end{eqnarray}
}\end{Def}
%
%
%
%%%%%%%%%%%%%%%%%%%%%%%%%%%%%%%%%%%%%%%%%%%%%%%%%%%%%%%%%%%%%%
%                 HORIZONTAL JACOBIAN
%%%%%%%%%%%%%%%%%%%%%%%%%%%%%%%%%%%%%%%%%%%%%%%%%%%%%%%%%%%%%%
%
%
%
\begin{Def}[Horizontal jacobian]\label{defhorjac}{\rm
Let $\Omega$ be an open subset of $\H^{n}$ and let $x\in\Omega$.  
The {\em horizontal jacobian} of a $C^1$ mapping
$f:\Omega\lra\R^k$ at $x$ is given by
\begin{eqnarray}
J_Hf(x)=\|\Lb_k\lls df(x)_{|H_x\H^{n}}\rls\|\,,
\end{eqnarray}
where $\Lb_k\lls df(x)_{|H_x\H^{n}}\rls:\Lb_k(H_x\H^{n})\lra
\Lb_k(\R^k)$.
}\end{Def}
%
%
%
%%%%%%%%%%%%%%%%%%%%%%%%%%%%%%%%%%%%%%%%%%%%%%%%%%%%%%%%%%%%%%
%           REMARK ON DEFINITION OF HORIZONTAL JACOBIAN
%%%%%%%%%%%%%%%%%%%%%%%%%%%%%%%%%%%%%%%%%%%%%%%%%%%%%%%%%%%%%%
%
%
%
From definition of horizontal jacobian, it follows that
it only depends on the restriction of $g$ to the horizontal subbundle, namely, from the ``sub-Riemannian metric".
Let us consider a horizontal frame $(Y_1,Y_2,\ldots,Y_{2n})$,
hence $J_Hf(x)$ is given by the jacobian of 
\begin{eqnarray}\label{hdm}
\nabla_Yf(x)=\left[\begin{array}{ccccccc}
Y_{1}f^{1}(x) & Y_{2}f^{1}(x) & \ldots & Y_{2n}f^{1}(x)  \\
Y_{1}f^{2}(x) & Y_{2}f^{2}(x) & \ldots & Y_{2n}f^{2}(x)  \\
\vdots        & \vdots        & \vdots & \vdots          \\
Y_{1}f^{k}(x) & Y_{2}f^{k}(x) & \ldots & Y_{2n}f^{k}(x)  \\
\end{array}\right].
\end{eqnarray}
As a consequence, we have the formula
\begin{eqnarray}\label{hjac}
J_Hf(x)=\sqrt{\sum_{1\leq i_{1}<i_{2}\cdots<i_{k}\leq 2n}
\left[M_{i_{1}i_{2}\cdots i_{k}}\left(\nabla_Yf(x)\right)\right]^{2}}.
\end{eqnarray}
%
%
%
%
%
%%%%%%%%%%%%%%%%%%%%%%%%%%%%%%%%%%%%%%%%%%%%%%%%%%%%%%%%%%%%%%
%            LEBESGUE MEASURE AND RIEMANNIAN VOLUME
%%%%%%%%%%%%%%%%%%%%%%%%%%%%%%%%%%%%%%%%%%%%%%%%%%%%%%%%%%%%%%
%
%
%
\begin{Pro}\label{egriem}
Let $(Y_1,Y_2,\ldots,Y_{2n},W)$ be an orthonormal frame
with respect to a left invariant metric $h$ and
let $F:\R^{2n+1}\lra\H^n$ define coordinates with respect
to this frame.
Then we have $F_\sharp\cL^{2n+1}=\vol_{2n+1}$,
where $\vol_{2n+1}$ denotes the Riemannian volume
measure with respect to the metric $h$.
\end{Pro}
{\sc Proof}.
Let $A$ be a measurable set of $\R^{2n+1}$. By classical
area formula and taking into account the left invariance
of both $\vol_p$ and $F_\sharp\cL^{2n+1}$ we have
$$
c\;\cL^{2n+1}(A)=\vol_{2n+1}(F(A))=\int_A J_hF(x)\,dx
$$
for some constant $c>0$. Then $\avint_AJ_hF=c$
for any measurable $A$. By continuity of $x\lra J_hF(x)$
we obtain that $J_hF(x)=c$ for any $x\in\R^q$.
We have $F=\exp\circ L$, with
\begin{eqnarray*}
L(y)=y_{2n+1}W+\sum_{i=1}^{2n}y_j\,Y_j	
\end{eqnarray*}
and $(Y_1,Y_2,\ldots,Y_{2n},W)$ is orthonormal.
Since the map $dF(0)=d\exp(0)\cir L=L$ has jacobian equal to one,
then $c=1$ and the thesis follows. $\Box$
\begin{Rem}{\rm
By previous proposition, the volume measure of a measurable subset $A$ of $\H^{n}$ corresponds to the 
$(2n+1)$-dimensional Lebesgue measure of the same subset
read with respect to coordinates associated to an orthonormal frame. Here the volume measure is defined by the same left invariant metric.
}\end{Rem}
Recall that $(F_\sharp\mu)(A)=\mu\lls F^{-1}(A)\rls$, where $\mu$
is a measure defined on the domain of $F$ and $A$ is a measurable
set defined on the codomain of $F$.
The $d$-dimensional spherical Hausdorff measure $\cS^{d}$ 
is defined as
\begin{eqnarray}
\cS^d(A)=\sup_{\ep>0}\;\inf\left\{\sum_{j=1}^\infty \diam(E_j)^d
\mid A\subset\bigcup E_j, \diam(E_j)\leq\ep\right\},
\end{eqnarray}
where the diameter is considered with respect to a
homogeneous distance $\rho$ of $\H^{n}$ and we do not consider any
dimensional factor.
The $k$-dimensional Hausdorff measure built with respect to the Riemannian distance is denoted by $\vol_k$ and it corresponds to the classical
Riemannian volume measure with respect to the graded metric $g$,
see for instance 3.2.46 of \cite{Fed}.
%
%
%
%%%%%%%%%%%%%%%%%%%%%%%%%%%%%%%%%%%%%%%%%%%%%%%%%%%%%%%%%%%%%%
%               HORIZONTAL P-VECTORS
%%%%%%%%%%%%%%%%%%%%%%%%%%%%%%%%%%%%%%%%%%%%%%%%%%%%%%%%%%%%%%
%
%
%
\begin{Def}[Horizontal $p$-vectors]{\rm
For each $x\in\H^{n}$, we say that any linear combination of wedge pro\-ducts $X_{j_1}(x)\wedge X_{j_2}(x)\wedge\cdots X_{j_p}(x)$, where $1\leq j_s\leq 2n$ and $j=1,\ldots,2n$, is
a {\em horizontal $p$-vector}. The space of horizontal $p$-vectors
is denoted by $\Lb_p(H_x\H^{n})$.
}\end{Def}
%
%
%
%%%%%%%%%%%%%%%%%%%%%%%%%%%%%%%%%%%%%%%%%%%%%%%%%%%%%%%%%%%%%%
%                 VERTICAL P-VECTORS
%%%%%%%%%%%%%%%%%%%%%%%%%%%%%%%%%%%%%%%%%%%%%%%%%%%%%%%%%%%%%%
%
%
%
\begin{Def}[Vertical $p$-vectors]{\rm
For each $x\in\H^{n}$, we say that any linear combination of wedge pro\-ducts $X_{j_1}(x)\wedge X_{j_2}(x)\wedge\cdots X_{j_{p-1}}(x)\wedge Z(x)$, where  $1\leq j_s\leq 2n$ and $j=1,\ldots,2n$, is a {\em vertical $p$-vector}. 
The space of vertical $p$-vectors is denoted by $\cV_p(H_x\H^{n})$.
}\end{Def}
%
%
%
%%%%%%%%%%%%%%%%%%%%%%%%%%%%%%%%%%%%%%%%%%%%%%%%%%%%%%%%%%%%%%
%            COMMENT ON VERTICAL P-VECTORS
%%%%%%%%%%%%%%%%%%%%%%%%%%%%%%%%%%%%%%%%%%%%%%%%%%%%%%%%%%%%%%
%
%
%
For every couple of simple $p$-vectors
$v_1\wedge\cdots\wedge v_p$, $w_1\wedge\cdots\wedge w_p\in \Lambda_p(T_x\H^n)$,
we define the scalar product induced by the left invariant
Riemannian metric $g$ on $T_x\H^n$ as
\begin{eqnarray}\label{pscalarproduct}
\lan v_1\wedge\cdots\wedge v_p,w_1\wedge\cdots\wedge w_p\ran=
\det\left(\lls g(x)(v_i,w_j\rls\right),
\end{eqnarray}
see for instance 1.7.5 of \cite{Fed} for more details.
This allows us to regard the space of vertical $p$-vectors
$\cV_p(T_x\H^{n})$ as the orthogonal complement of the horizontal subspace $\Lb_p(H_x\H^{n})$.
We have the orthogonal decomposition
\begin{eqnarray}\label{orthdec}
\Lb_p(T_x\H^{n})=\Lb_p(H_x\H^{n})\oplus\cV_p(T_x\H^{n}),
\end{eqnarray}
which generalizes the case $p=1$, corresponding to 
$T_x\H^{n}=H_x\H^{n}\oplus\lan Z(x)\ran$.
%
%
%
%%%%%%%%%%%%%%%%%%%%%%%%%%%%%%%%%%%%%%%%%%%%%%%%%%%%%%%%%%%%%%
%                VERTICAL PROJECTION
%%%%%%%%%%%%%%%%%%%%%%%%%%%%%%%%%%%%%%%%%%%%%%%%%%%%%%%%%%%%%%
%
%
%
\begin{Def}[Vertical projection]{\rm
Let $x\in\H^{n}$ and let $\xi\in\Lb_p(T_x\H^{n})$.
The orthogonal decomposition $\xi=\xi_H+\xi_\cV$ associated
to (\ref{orthdec}) uniquely defines the vertical $p$-vector $\xi_\cV\in\cV_p(T_x\H^{n})$.
We say that $\xi_\cV$ is the {\em vertical projection} of $\xi$
and that the mapping $\pi_\cV:\Lb_p(T_x\H^{n})\lra\cV_p(T_x\H^{n})$,
which associates $\xi_\cV$ to  $\xi$, is 
the {\em vertical projection}.
}\end{Def}
%
%
%
%%%%%%%%%%%%%%%%%%%%%%%%%%%%%%%%%%%%%%%%%%%%%%%%%%%%%%%%%%%%%%
%           COMMENT ON THE VERTICAL PROJECTION
%%%%%%%%%%%%%%%%%%%%%%%%%%%%%%%%%%%%%%%%%%%%%%%%%%%%%%%%%%%%%%
%
%
%
We have omitted $x$ in the definition
of vertical projection $\pi_\cV$.
%
%
%
%%%%%%%%%%%%%%%%%%%%%%%%%%%%%%%%%%%%%%%%%%%%%%%%%%%%%%%%%%%%%%
%                    CHARACTERISTIC POINT
%%%%%%%%%%%%%%%%%%%%%%%%%%%%%%%%%%%%%%%%%%%%%%%%%%%%%%%%%%%%%%
%
%
%
\begin{Def}[Characteristic points and transverse points]\label{defchar}
{\rm Let $\Sigma\subset\Omega$ be a $C^1$ submanifold and let $x\in\Sigma$.
We say that $x\in\Sigma$ is a {\em characteristic point} if $T_x\Sigma\subset H_x\H^{n}$ and that it is a {\em transverse point} otherwise.
The {\em characteristic set} of $\Sigma$ is the subset of all characteristic points and it is denoted by $C(\Sigma)$.
}\end{Def}
%
%
%
%%%%%%%%%%%%%%%%%%%%%%%%%%%%%%%%%%%%%%%%%%%%%%%%%%%%%%%%%%%%%%
%                TANGENT P-VECTOR 
%%%%%%%%%%%%%%%%%%%%%%%%%%%%%%%%%%%%%%%%%%%%%%%%%%%%%%%%%%%%%%
%
%
%
Recall that a {\em tangent $p$-vector} to a $p$-dimensional
submanifold $\Sigma$ of class $C^1$ at $x\in\Sigma$ is
defined by the wedge product
$t_1\wedge t_2\wedge\cdots\wedge t_p$, where
$(t_1,\ldots,t_p)$ is an orthonormal basis of $T_x\Sigma$.
We denote this simple $p$-vector by $\tau_\Sigma(x)$.
%
%
%
%%%%%%%%%%%%%%%%%%%%%%%%%%%%%%%%%%%%%%%%%%%%%%%%%%%%%%%%%%%%%%
%                REMARK ON TANGENT P-VECTORS
%%%%%%%%%%%%%%%%%%%%%%%%%%%%%%%%%%%%%%%%%%%%%%%%%%%%%%%%%%%%%%
%
%
%
Notice that the tangent $p$-vector (which belongs to a one-dimensional space) cannot be continuously defined on all of $\Sigma$,
unless the submanifold is oriented.
%
%
%
%%%%%%%%%%%%%%%%%%%%%%%%%%%%%%%%%%%%%%%%%%%%%%%%%%%%%%%%%%%%%%
%             VERTICAL TANGENT P-VECTOR %%%%%%%%%%%%%%%%%%%%%%%%%%%%%%%%%%%%%%%%%%%%%%%%%%%%%%%%%%%%%%
%
%
%
\begin{Def}[Vertical tangent $p$-vector]\label{vertangpvect}
{\rm Let $\Sigma\subset\Omega$ be a $p$-dimensional submanifold of class $C^1$ and let $x\in\Sigma$.
A {\em vertical tangent $p$-vector} to $\Sigma$ at $x$ is 
defined by $\pi_\cV(\tau_\Sigma)$, where $\tau_\Sigma$ is a
tangent $p$-vector and $\pi_\cV$ is the vertical projection.
The vertical tangent $p$-vector will be denoted by $\tau_{\Sigma,\cV}(x)$.
}\end{Def}
%%%%%%%%%%%%%%%%%%%%%%%%%%%%%%%%%%%%%%%%%%%%%%%%%%%%%%%%%%%%%%
%%%%%%%%%%%%%%%%%%%%%%%%%%%%%%%%%%%%%%%%%%%%%%%%%%%%%%%%%%%%%%
%%
%%
%%
%%
%%
\section{Blow-up at transverse points}\label{blwtran}
%%
%%
%%
%%
%%
%%%%%%%%%%%%%%%%%%%%%%%%%%%%%%%%%%%%%%%%%%%%%%%%%%%%%%%%%%%%%%
%%%%%%%%%%%%%%%%%%%%%%%%%%%%%%%%%%%%%%%%%%%%%%%%%%%%%%%%%%%%%%
%
%
% 
This section is devoted to the proof of Theorem~\ref{bwpt}.
In the following proposition, we give a simple characterization of characteristic points using vertical tangent $p$-vectors.
%
%
%
%%%%%%%%%%%%%%%%%%%%%%%%%%%%%%%%%%%%%%%%%%%%%%%%%%%%%%%%%%%%%%
%          CHARATERIZATION OF CHARACTERISTIC POINTS
%          THROUGH VERTICAL PROJECTION
%%%%%%%%%%%%%%%%%%%%%%%%%%%%%%%%%%%%%%%%%%%%%%%%%%%%%%%%%%%%%%
%
%
%
\begin{Pro}\label{charcharvp}
Let $\Sigma\subset\Omega$ be a submanifold of class $C^1$
and let $x\in\Sigma$. Then $x\in C(\Sigma)$ if and only if
$\tau_{\Sigma,\cV}(x)=0$.
\end{Pro}
{\sc Proof}.
Let $x\in\Sigma$ and let $(t_1,t_2,\ldots,t_p)$ be an orthonormal
basis of $T_x\Sigma$. We have the unique decomposition
$
t_j=V_j+\gm_j\;Z,
$
where $V_j\in H_x\H^{n}$ for every $j=1,\ldots,p$.
It follows that
\begin{eqnarray*}
&&\tau=t_1\wedge t_2\wedge\cdots\wedge t_p\\
&&=\lls V_1+\gm_1 Z\rls
\wedge\lls V_2+\gm_2 Z\rls\wedge\cdots
\wedge\lls V_p+\gm_p Z\rls\\
&&=V_1\wedge V_2\wedge\cdots\wedge V_p+\sum_{j=1}^{p}\gm_j\;\;
V_1\wedge V_2\wedge\cdots V_{j-1}\wedge Z\wedge V_{j+1}
\wedge\cdots\wedge V_{p}.
\end{eqnarray*}
Assume that $x\notin C(\Sigma)$. If $V_1,V_2,\ldots,V_p$ are linearly dependent, then we get
\begin{eqnarray*}
&&t_1\wedge t_2\wedge\cdots\wedge t_p
=\sum_{j=1}^{p}\gm_j\;\;
V_1\wedge V_2\wedge\cdots V_{j-1}\wedge Z\wedge V_{j+1}
\wedge\cdots\wedge V_{p}.
\end{eqnarray*}
As a result, $\pi_\cV(\tau)=\tau$ hence it is not vanishing.
If $V_1,\ldots,V_p$ are linearly independent, then
all wedge products of the form
\begin{eqnarray}
V_1\wedge V_2\wedge\cdots V_{j-1}\wedge Z\wedge V_{j+1}\wedge\cdots\wedge V_p
\end{eqnarray}
are non-vanishing for every $j=1,\ldots,p$.
The fact that $x$ is transverse implies that there exists
$\gm_{j_0}\neq0$, then the projection
\begin{eqnarray}\label{projv}
\pi_\cV(\tau)=\sum_{j=1}^p\,\gm_j\;	V_1\wedge V_2\wedge\cdots V_{j-1}\wedge Z\wedge V_{j+1}\wedge\cdots\wedge V_p
\end{eqnarray}
is non-vanishing. Conversely, if $\pi_\cV(\tau)\neq0$, then (\ref{projv})
yields some $\gm_{j_1}\neq0$, therefore
$t_{j_1}\notin H_x\H^n$. $\Box$
%
%
%
%%%%%%%%%%%%%%%%%%%%%%%%%%%%%%%%%%%%%%%%%%%%%%%%%%%%%%%%%%%%%%
%          CHARATERIZATION OF CHARACTERISTIC POINTS
%          THROUGH HORIZONTAL JACOBIAN
%%%%%%%%%%%%%%%%%%%%%%%%%%%%%%%%%%%%%%%%%%%%%%%%%%%%%%%%%%%%%%
%
%
%
\begin{Pro}\label{charcharhj}
Let $f:\Omega\lra\R^k$ be of class $C^1$, with surjective differential at each point of $\Omega$. Let $\Sigma$ denote the submanifold $f^{-1}(0)$ of $\Omega$ and let $x\in\Sigma$.
Then $x\in C(\Sigma)$ if and only if
$df(x)_{|H_x\H^{n}}$ is not surjective.
\end{Pro}
{\sc Proof}.
We first notice that 
$
\ker\Big(df(x)_{|H_x\H^{n}}\Big)=T_x\Sigma\cap H_x\H^{n},
$ 
then we have
\begin{eqnarray}\label{kkkeyeq}
\dim(H_x\H^{n}\cap T_x\Sigma)=2n
-\dim\Big(\im\lls df(x)_{|H_x\H^{n}}\rls\Big).
\end{eqnarray}
This last formula allows us to get our claim as follows.
Assume that $x\in C(\Sigma)$. Then $T_x\Sigma\subset H_x\H^{n}$ and (\ref{kkkeyeq}) gives
$$
2n+1-k=2n-\dim\Big(\im\lls df(x)_{|H_x\H^{n}}\rls\Big).
$$
From this equation we conclude that $df(x)_{|H_x\H^{n}}$
is not surjective.
Conversely, if $df(x)_{|H_x\H^{n}}$ is not surjective, then (\ref{kkkeyeq})
implies 
\begin{eqnarray*}
\dim(H_x\H^{n}\cap T_x\Sigma)\geq 2n-k+1=\dim(T_x\Sigma)
\end{eqnarray*}
therefore $T_x\Sigma\subset H_x\H^{n}$, namely, $x\in C(\Sigma)$. $\Box$

%
%
%
%%%%%%%%%%%%%%%%%%%%%%%%%%%%%%%%%%%%%%%%%%%%%%%%%%%%%%%%%%%%%%
%                   QUOTIENT OF JACOBIANS
%%%%%%%%%%%%%%%%%%%%%%%%%%%%%%%%%%%%%%%%%%%%%%%%%%%%%%%%%%%%%%
%
%
%
\begin{The}\label{jhjg}
Let $f:\Omega\lra\R^k$ be of class $C^1$, with surjective differential
at each point of $\Omega$. Let $\Sigma$ denote the submanifold
$f^{-1}(0)$ of $\Omega$ and let $x\in\Sigma$.
Then we have
\begin{eqnarray}\label{jhoverjg}
|\tau_{\Sigma,\cV}(x)|=\frac{J_Hf(x)}{J_gf(x)}\;.
\end{eqnarray}
\end{The}
{\sc Proof}.
Left invariance of Riemannian metric allows us to
consider the left translated submanifold $l_{x^{-1}}\Sigma$.
Replacing $f$ with $f\circ l_x$ and $\Omega$ with $l_{x^{-1}}\Omega$
we can assume that $x$ is the unit element $0$ of $\H^{n}$.
Recall that $l_x:\H^{n}\lra\H^{n}$ is the left translation
$l_x(y)=x\cdot y$.
If $x\in C(\Sigma)$, then Proposition~\ref{charcharvp} and Proposition~\ref{charcharhj} make (\ref{jhoverjg}) the trivial
identity $0=0$. 
Assume that $x\in\Sigma\sm C(\Sigma)$.
Then Proposition~\ref{charcharhj} implies that the horizotal gradients
$$
\nabla_Hf^i=(X_1f^i(0),X_2f^i(0),\ldots,X_{2n}f^i(0))
\quad\mbox{for}\quad i=1,2,\ldots,k
$$
span a $k$-dimensional space of $\R^{2n}$.
Let $c_1,c_2,\ldots,c_k\in\R^{2n}$ be orthogonal unit vectors
generating this vector space and choose
$c_{k+1},\ldots,c_{2n}\in\R^{2n}$ such that $(c_1,c_2,\ldots,c_{2n})$
is an orthonormal basis of $\R^{2n}$. 
These vectors allow us to define a new horizontal frame
\begin{eqnarray}\label{keyrotation1}
Y_j=\sum_{k=1}^{2n}c_j^k\; X_k\quad\mbox{for every}\quad j=1,\ldots,2n.
\end{eqnarray}
We denote by $C$ the $2n\times2n$ orthogonal matrix whose $i$-th column
corresponds to the vector $c_i$, then by our choice of vectors $c_j$, we obtain $\nabla_Yf(x)=\nabla_Xf(x)\;C$ and 
\begin{eqnarray}\label{difhmtx1}
\nabla_Yf(x)C=
\left[\begin{array}{ccccccc}
\lan \nabla_H f^1,c_1\ran &\lan \nabla_H f^1,c_2\ran & \cdots 
&\lan \nabla_H f^1,c_k\ran & 0 & \cdots & 0\\ 
\lan \nabla_H f^2,c_1\ran &\lan \nabla_H f^2,c_2\ran & \cdots 
&\lan \nabla_H f^2,c_k\ran & 0 & \cdots & 0\\ 
\vdots  & \vdots & \cdots & \vdots & 0 & \cdots & 0\\ 
\lan \nabla_H f^k,c_1\ran &\lan \nabla_H f^k,c_2\ran & \cdots 
&\lan \nabla_H f^k,c_k\ran & 0 & \cdots & 0\\ 
\end{array}\right],
\end{eqnarray}
where the symbol $\lan,\ran$ denotes the standard scalar product
of $\R^{2n}$. Let us consider $F:\R^{2n+1}\lra\H^{n}$, defining graded coordinates $(y_1,\ldots,y_{2n+1})$ associated to the frame $(Y_1,\ldots,Y_{2n},Z)$,
according to Definition~\ref{gradef}.
Then the differential of $f$ at $0$ with respect to
$(y_1,\ldots,y_{2n+1})$ can be represented by the matrix
\begin{eqnarray}\label{dforig1}
\nabla_yf(0)=\left[\begin{array}{cccccccc}
f^1_{y_1}(0) & f^1_{y_2}(0) & \cdots & f^1_{y_k}(0) & 0 & \cdots & 0 &
f^1_{y_{2n+1}}(0)\\ 
f^2_{y_1}(0) & f^2_{y_2}(0) & \cdots & f^1_{y_k}(0) & 0 & \cdots & 0 &
f^2_{y_{2n+1}}(0)\\ 
\vdots  & \vdots & \cdots & \vdots & 0 & \cdots & 0 & \vdots\\ 
f^k_{y_1}(0) & f^k_{y_2}(0) & \cdots & f^k_{y_k}(0) & 0 & \cdots & 0 &
f^k_{y_{2n+1}}(0) \\ 
\end{array}\right].
\end{eqnarray}
It follows that 
$$
f^i_{y_j}(0)=\lan \nabla_H f^i,c_j\ran
$$
for every $i,j=1,\ldots,2n$.
The implicit function theorem gives us
a $C^1$ map $\ph:A\lra\R^k$ such that $A\subset\R^p$ is an
open neighbourhood of the origin and 
\begin{eqnarray}\label{impeq1}
f(\ph^1(\tilde{y}),\ldots,\ph^k(\tilde{y}),y_{k+1},\ldots,y_{2n+1})=0
\end{eqnarray}
for every $\tilde{y}=(y_{k+1},\ldots,y_{2n+1})\in A$.
Then we define the mapping $\phi:A\lra\R^{2n+1}$ as
\begin{eqnarray}\label{phiform}
\phi(\tilde{y})=(\ph^1(\tilde{y}),\ldots,\ph^k(\tilde{y}),y_{k+1},\ldots,y_{2n+1}),
\end{eqnarray}
so that differentiating (\ref{impeq1}) we get
\begin{eqnarray}\label{0dfimp1}
0=\der_{y_j}(f^i\circ\phi)
=\sum_{l=1}^kf^i_{y_l}\,\ph^l_{y_j}+f^i_{y_j}
\end{eqnarray}
for every $i=1,\ldots,k$ and $j=k+1,\ldots,2n+1$.
Equations (\ref{0dfimp1}) can be more concisely written
in matrix form as follows
\begin{eqnarray}\label{matrixform1}
\nabla_zf\;\;\ph_{y_j}=-f_{y_j},
\end{eqnarray}
where $z=(y_1,\ldots,y_k)$, the $k\times k$ matrix $\nabla_zf$ has
coefficients $f^i_{y_l}$, where $i,l=1,\ldots,k$
and $j=k+1,\ldots,2n+1$.
In order to achieve a more explicit formula for the differential of
the implicit map, we explicitly write the inverse matrix of 
$\nabla_zf$ as 
\begin{eqnarray*}
\left(\nabla_zf\right)^{-1}=
\frac{1}{M_{12\cdots k}(\nabla_zf)}
\left[\begin{array}{cccc}
C_{11}(\nabla_zf) & C_{21}(\nabla_zf) & \cdots & C_{k1}(\nabla_zf) \\
C_{12}(\nabla_zf) & C_{22}(\nabla_zf) & \cdots & C_{k1}(\nabla_zf) \\
\vdots     & \vdots    &  \cdots & \vdots \\
C_{1k}(\nabla_zf) & C_{2k}(\nabla_zf) & \cdots & C_{kk}(\nabla_zf) \\
\end{array}\right]\;,
\end{eqnarray*}
where $C_{ij}(\nabla_zf)$ denotes the cofactor of $\nabla_zf$,
which is equal to $(-1)^{i+j}\det(\mbox{{\rm \^D}}_{ij}f)$ and
$\mbox{{\em \^D}}_{ij}f$ is the $(k-1)\times(k-1)$ square matrix obtained
by removing the $i$-th row and the $j$-th column from $\nabla_zf$.
In view of (\ref{matrixform1}) we have
\begin{eqnarray*}%\label{}
\ph_{y_j}=-\left(\nabla_zf\right)^{-1}f_{y_j}
=-\frac{1}{M_{12\cdots k}(\nabla_zf)}
\left[\begin{array}{c}
\sum_{i=1}^k C_{i1}(\nabla_zf)f^i_{y_j} \\
\sum_{i=1}^k C_{i2}(\nabla_zf)f^i_{y_j} \\
\vdots    \\
\sum_{i=1}^k C_{ik}(\nabla_zf)f^i_{y_j}
\end{array}\right].
\end{eqnarray*}
An elementary formula for computing the determinant of a matrix implies
\begin{eqnarray*}
\sum_{i=1}^k C_{is}(\nabla_zf)f^i_{y_j}=M_{12\cdots s-1\,j\,s+1\cdots k}(\nabla_yf)
\end{eqnarray*}
for every $j=k+1,\ldots,2n+1$.
As a consequence, we get
\begin{eqnarray}\label{algform1}
\ph^s_{y_j}=-\frac{M_{12\cdots s-1\,j\,s+1\cdots k}(\nabla_yf)}
{M_{1\cdots k}(\nabla_zf)}.
\end{eqnarray}
Note that $M_{1\cdots k}(\nabla_zf)$ corresponds to the determinant
of the matrix $\nabla_zf$. As a consequece of (\ref{algform1})
and of (\ref{dforig1}), we conclude that
\begin{eqnarray*}%\label{}
\ph^s_{y_j}(0)=0
\end{eqnarray*}
for every $j=k+1,\ldots,2n$. Previous considerations and expression
(\ref{phiform}) lead us to the formula
\begin{eqnarray}\label{dphi01}
\nabla_{\tilde{y}}\phi(0)
=\left[\begin{array}{ccccc}
0\qquad      &  0\qquad     & \cdots\qquad & 0\,       & \ph^1_{y_{2n+1}}(0) \\
0\qquad      &  0\qquad     & \cdots\qquad & 0\,       & \ph^2_{y_{2n+1}}(0) \\
\vdots\qquad & \vdots\qquad & \ddots\qquad & \vdots\,  & \vdots      \\
0\qquad      &  0\qquad     & \cdots\qquad & 0\,       & \ph^k_{y_{2n+1}}(0) \\
1\qquad      &  0\qquad     & \cdots\qquad & 0\,       &   0      \\
0\qquad      &  1\qquad     & \cdots\qquad & 0\,       &   0 \\
\vdots\qquad &  0\qquad     & \ddots\qquad & 0\,       &   0  \\
\vdots\qquad &  \vdots\qquad& \ddots\qquad & 1\,       & \vdots \\
0\qquad      &  0\qquad     & \cdots\qquad & 0\,       & 1
\end{array}\right],
\end{eqnarray}
where $\nabla_{\ytil}\phi(0)$ is a $(2n+1)\times p$ matrix 
whose $p\times p$ lower block is the identity matrix.
Notice that columns of (\ref{dphi01}) represent a basis of
the tangent space $T_0\Sigma$ with respect to coordinates
$(y_{k+1},\ldots,y_{2n+1})$. More precisely, the set of vectors
\begin{eqnarray*}%\label{}
\left(Y_{k+1}(0),Y_{k+2}(0),\ldots,Y_{2n}(0),\frac{Z(0)+\sum_{j=1}^k v_j\,Y_j(0)}{\Big(1+\sum_{j=1}^k\;v_j^2\Big)^{1/2}}\right)
\end{eqnarray*}
form an orthonormal basis of $T_0\Sigma$, where we have defined
$v_j=\ph^j_{y_{2n+1}}(0)$.
Then the tangent $p$-vector $\tau_\Sigma$ to $\Sigma$ at $0$ is given by the wedge product
\begin{eqnarray*}%\label{}
\tau_{\Sigma}(0)=\frac{Y_{k+1}(0)\wedge Y_{k+2}(0)\wedge\cdots\wedge
Y_{2n}(0)\wedge\Big(Z(0)+\sum_{j=1}^k v_j\,Y_j(0)\Big)}{\Big(1+\sum_{j=1}^k\;v_j^2\Big)^{1/2}}\;.
\end{eqnarray*}
Obviously, $p$-vectors $Y_{k+1}(0)\wedge Y_{k+2}(0)\wedge\cdots\wedge
Y_{2n}(0)\wedge Y_j(0)$ are horizontal, hence they disappear
in the vertical projection. It follows that 
\begin{eqnarray*}%\label{}
\tau_{\Sigma,\cV}(0)=\pi_\cV\lls\tau_{\Sigma,\cV}(0)\rls=\frac{Y_{k+1}(0)\wedge Y_{k+2}(0)\wedge\cdots\wedge Y_{2n}(0)\wedge Z(0)}{\Big(1+\sum_{j=1}^k\;v_j^2\Big)^{1/2}}\;
\;,
\end{eqnarray*}
therefore we clearly obtain 
\begin{eqnarray}\label{verttang}
|\tau_{\Sigma,\cV}(0)|=\Big(1+\sum_{l=1}^k v_l^2\Big)^{-1/2}.
\end{eqnarray}
Due to formula (\ref{algform1}) in the case $j=2n+1$ and to (\ref{dforig1}), we obtain
\begin{eqnarray*}
1+\sum_{l=1}^k v_l^2=\frac{\lls M_{1\cdots k}(\nabla_zf)\rls^2+\sum_{s=1}^k\lls M_{12\cdots s-1\,2n+1\,s+1\cdots k}
(\nabla_yf)\rls^2}{\lls M_{1\cdots k}(\nabla_zf)\rls^2}
=\bigg(\frac{J_gf(0)}{J_Hf(0)}\bigg)^2\,,
\end{eqnarray*}
then (\ref{verttang}) shows the valdity of (\ref{jhoverjg}) in the
case $x=0$. Left invariance of the Riemannian metric $g$
leads us to the conclusion. $\Box$
%
%
%
%%%%%%%%%%%%%%%%%%%%%%%%%%%%%%%%%%%%%%%%%%%%%%%%%%%%%%%%%%%%%
%                    METRIC FACTOR
%%%%%%%%%%%%%%%%%%%%%%%%%%%%%%%%%%%%%%%%%%%%%%%%%%%%%%%%%%%%%
%
%
%
\begin{Def}[Metric factor]\label{metrfct}{\rm
Let $\tau$ be a vertical simple $p$-vector of $\Lb_p(\ch^{n})$ and let
$\cL(\tau)$ be the unique associated subspace, with $L=\exp\cL(\tau)$.
The {\em metric factor} of a homogeneous distance $\rho$ with respect to
$\tau$ is defined by
\begin{eqnarray*}
\theta_p^\rho(\tau)=\mcH^p_{|\cdot|}\left(F^{-1}(L\cap B_1)\right),
\end{eqnarray*}
where $F:\R^{2n+1}\lra\H^{n}$ defines a system of graded coordinates,
$\mcH^p_{|\cdot|}$ denotes the $p$-dimensional Hausdorff measure with respect to the Euclidean distance of $\R^{2n+1}$ and $B_1$ is the unit ball of 
$\H^{n}$ with respect to the distance $\rho$.
Recall that the subspace associated to a simple $p$-vector $\tau$
is defined as $\{v\in\ch^{n}\mid v\wedge\tau=0 \}$.
}\end{Def}
\begin{Rem}{\rm
In the case of subspaces $\cL$ of codimension one, the notion of metric
factor fits into the one introduced in \cite{Mag2}. 
It is easy to observe that the notion of metric factor does not depend on the system of coordinates we are using. In fact,
$F_1^{-1}\circ F_2:\R^{2n+1}\lra\R^{2n+1}$ is an Euclidean isometry
whenever $F_1,F_2:\R^{2n+1}\lra\H^{n}$ represent systems of 
graded coordinates with respect to the same left invariant Riemannian metric.
}\end{Rem}
%
%
%
%%%%%%%%%%%%%%%%%%%%%%%%%%%%%%%%%%%%%%%%%%%%%%%%%%%%%%%%%%%%%
%         PROOF OF BLOW-UP OF CODIMENSION k SUBMANIFOLDS
%%%%%%%%%%%%%%%%%%%%%%%%%%%%%%%%%%%%%%%%%%%%%%%%%%%%%%%%%%%%%
%
%
%
{\em Proof of Theorem~\ref{bwpt}}.
As in the proof of Theorem~\ref{jhjg}, left invariance of the   
Riemannian metric $g$ allows us to assume that $x=0$.
For $r_0>0$ sufficiently small, we can suppose the existence
of a function $f:B_{r_0}\lra\R^k$ such that
$\Sigma\cap B_{r_0}=f^{-1}(0)$ and whose differential is surjective
at every point of $B_{r_0}$.
By Proposition~\ref{charcharhj}, the horizontal gradients
$$
\nabla_Hf^i=(X_1f^i(0),X_2f^i(0),\ldots,X_{2n}f^i(0))
\quad\mbox{for}\quad i=1,2,\ldots,k
$$
span a $k$-dimensional space of $\R^{2n}$.
Now, repeating the argument in the proof of Theorem~\ref{jhjg},
we define the system of graded coordinates $(y_1,\ldots,y_{2n+1})$
associated to the frame $(Y_1,\ldots,Y_{2n},Z)$, 
where $Y_j$ are given by (\ref{keyrotation1}).
The differential of $f$ at $0$ can be represented by the matrix
\begin{eqnarray}\label{dforig}
\nabla_yf(0)=\left[\begin{array}{cccccccc}
f^1_{y_1}(0) & f^1_{y_2}(0) & \cdots & f^1_{y_k}(0) & 0 & \cdots & 0 &
f^1_{y_{2n+1}}(0)\\ 
f^2_{y_1}(0) & f^2_{y_2}(0) & \cdots & f^1_{y_k}(0) & 0 & \cdots & 0 &
f^2_{y_{2n+1}}(0)\\ 
\vdots  & \vdots & \cdots & \vdots & 0 & \cdots & 0 & \vdots\\ 
f^k_{y_1}(0) & f^k_{y_2}(0) & \cdots & f^k_{y_k}(0) & 0 & \cdots & 0 &
f^k_{y_{2n+1}}(0) \\ 
\end{array}\right]\,,
\end{eqnarray}
whose first $k$ columns are linearly independent.
By the implicit function theorem ther eexists 
a $C^1$ mapping $\ph:A\lra\R^k$ such that $A\subset\R^p$ is an
open neighbourhood of the origin and 
\begin{eqnarray}\label{impeq}
f(\ph^1(\tilde{y}),\ldots,\ph^k(\tilde{y}),y_{k+1},\ldots,y_{2n+1})=0
\end{eqnarray}
for every $\tilde{y}=(y_{k+1},\ldots,y_{2n+1})\in A$.
Proceeding as in the proof of Theorem~\ref{jhjg}, we define
the mapping $\phi:A\lra\R^{2n+1}$ as
\begin{eqnarray}%\label{}
\phi(\tilde{y})=(\ph^1(\tilde{y}),\ldots,\ph^k(\tilde{y}),y_{k+1},\ldots,y_{2n+1}),
\end{eqnarray}
and by the same computations, differentiating (\ref{impeq}) we obtain 
\begin{eqnarray}\label{dphi0}
\nabla_{\tilde{y}}\phi(0)
=\left[\begin{array}{ccccc}
0\qquad      &  0\qquad     & \cdots\qquad & 0\,       & \ph^1_{y_{2n+1}}(0) \\
0\qquad      &  0\qquad     & \cdots\qquad & 0\,       & \ph^2_{y_{2n+1}}(0) \\
\vdots\qquad & \vdots\qquad & \ddots\qquad & \vdots\,  & \vdots      \\
0\qquad      &  0\qquad     & \cdots\qquad & 0\,       & \ph^k_{y_{2n+1}}(0) \\
1\qquad      &  0\qquad     & \cdots\qquad & 0\,       &   0      \\
0\qquad      &  1\qquad     & \cdots\qquad & 0\,       &   0 \\
\vdots\qquad &  0\qquad     & \ddots\qquad & 0\,       &   0  \\
\vdots\qquad &  \vdots\qquad& \ddots\qquad & 1\,       & \vdots \\
0\qquad      &  0\qquad     & \cdots\qquad & 0\,       & 1
\end{array}\right],
\end{eqnarray}
where $\nabla_{\ytil}\phi(0)$ is a $(2n+1)\times p$ matrix 
whose $p\times p$ lower block is the identity matrix.
For each $r<r_0$, write the ball $B_r$ in terms of graded coordinates defining $\Btil_r=F^{-1}(B_r)\subset\R^{2n+1}$.
The surface $\Sigma$ read in graded coordinates can be seen
as the image of $\phi$.
Then we have established
\begin{eqnarray}\label{quotblw}
\frac{\vol_p(\Sigma\cap B_r)}{r^{p+1}}=r^{-1-p}
\int_{\phi^{-1}(\Btil_r)}J_g\phi(\ytil)\,d\ytil\,.
\end{eqnarray}
The dilation $\delta_r$ restricted to coordinates $(y_1,\ldots,y_{2n+1})$ gives
\begin{eqnarray}\label{expldil}
\delta_r\ytil=\delta_r\lls(y_{k+1},\ldots,y_{2n+1})\rls
=(ry_{k+1},ry_{k+2},\ldots,ry_{2n},r^2y_{2n+1}),
\end{eqnarray}
therefore, performing a change of variable in (\ref{quotblw}) we get
\begin{eqnarray}\label{quotbl}
\frac{\vol_p(\Sigma\cap B_{x,r})}{r^{p+1}}=
\int_{\delta_{1/r}\lls\phi^{-1}(\Btil_r)\rls}
J_g\phi(\delta_r\ytil)\,d\ytil\,.
\end{eqnarray}
The set $\delta_{1/r}\lls\phi^{-1}(\Btil_r)\rls$ can be written as follows
\begin{eqnarray}%\label{}
(\delta_{1/r}\circ\phi\circ\delta_r)^{-1}(\Btil_1)
=\left\{\ytil\in\R^p\bigg|
\left(\frac{\ph^1(\delta_r\ytil)}{r},\ldots,\frac{\ph^k(\delta_r\ytil)}{r},y_{k+1}
,\ldots,y_{2n+1}\right)\in\Btil_1\right\}.
\end{eqnarray}
From expressions (\ref{dphi0}) and  (\ref{expldil}) one easily gets that
\begin{eqnarray}%\label{}
\lim_{r\ra0^+}\frac{\ph^j(\delta_r\ytil)}{r}=0
\end{eqnarray}
for every $j=1,\ldots,k$. As a result, the limit
\begin{eqnarray}%\label{}
{\bf 1}_{\delta_{1/r}\lls\phi^{-1}(\Btil_r)\rls}\lra{\bf 1}_{\Btil_1\cap\Pi}
\quad\mbox{as}\quad r\ra0^+\end{eqnarray}
holds a.e. in $\R^p$, where we have defined 
\begin{eqnarray*}%\label{}
\Pi=\left\{(0,\ldots,0,y_{k+1},\ldots,y_{2n+1})\in\R^{2n+1}\mid y_j\in\R, j=k+1,\ldots,2n+1\right\}.
\end{eqnarray*}
From (\ref{quotbl}), we conclude that
\begin{eqnarray}\label{limitjactheta}
\lim_{r\ra0^+}\frac{\vol_p(\Sigma\cap B_{x,r})}{r^{p+1}}=J_g\phi(0)\; \mcH^p(\Pi\cap\Btil_1).
\end{eqnarray}
To compute $J_g\phi(0)$, we use both the canonical form of the tangent space $T_0\Sigma$ given by (\ref{dphi0}) and the fact that our frame $(Y_1,\ldots,Y_{2n},Z)$ is orthonormal.
Thus, according to Definition~\ref{riemjacobian} 
the Riemannian jacobian of $\phi$ at zero is given by
\begin{eqnarray}\label{jaconcan}
J_g\phi(0)=\bigg(1+\sum_{j=1}^k\;v_j^2\bigg)^{1/2},
\end{eqnarray}
where we have defined $v_j=\ph_{y_{2n+1}}^j(0)$ for every $j=1,\ldots,2n$. Again, following the same steps of the proof of Theorem~\ref{jhjg}, we get
\begin{eqnarray*}
|\tau_{\Sigma,\cV}(0)|=\Big(1+\sum_{l=1}^k v_l^2\Big)^{-1/2}=
\lls J_g\phi(0)\rls^{-1}.
\end{eqnarray*}
Then (\ref{limitjactheta}) yields
\begin{eqnarray}\label{almfin}
\lim_{r\ra0^+}\frac{\vol_p(\Sigma\cap B_r)}{r^{p+1}}=\frac{\mcH^p(\Pi\cap\Btil_1)}
{|\tau_{\Sigma,\cV}(0)|}.
\end{eqnarray}
The subspace $\cL\lls\tau_{\Sigma,\cV}(x)\rls$ associated to the $p$-vector 
$\tau_{\Sigma,\cV}(x)$ satisfies the relation
$$
\exp\left(\cL\lls\tau_{\Sigma,\cV}(x)\rls\right)=F(\Pi)
$$
therefore the metric factor of $\rho$ with respect to $\tau_{\Sigma,\cV}(x)$ is
$\mcH^p(\Pi\cap\Btil_1)$.
This fact along with (\ref{almfin}) implies the validity of (\ref{keyblow}) and ends the proof. $\Box$
%
%
%
%
%%%%%%%%%%%%%%%%%%%%%%%%%%%%%%%%%%%%%%%%%%%%%%%%%%%%%%%%%%%%%
%
%    SECTION ON BLOW-UP OF CODIMENSION k SUBMANIFOLDS
%
%%%%%%%%%%%%%%%%%%%%%%%%%%%%%%%%%%%%%%%%%%%%%%%%%%%%%%%%%%%%%
%
%
%
%
\section{Spherical Hausdorff measure of submanifolds}\label{sphersub}

This section deals with various applications of Theorem~\ref{constmetintsp}. 
A key result to obtain this theorem is the $\cS^{Q-k}$-negligibility
of characteristic points of a $k$-codimensional submanifold
of a Carnot group of Hausdorff dimension $Q$, see \cite{Mag5}.
This result in the case of Heisenberg groups reads as follows.
%
%
%
%%%%%%%%%%%%%%%%%%%%%%%%%%%%%%%%%%%%%%%%%%%%%%%%%%%%%%%%%%%%%%
%        NEGLIGIBILITY OF CHARACTERISTIC POINTS  %%%%%%%%%%%%%%%%%%%%%%%%%%%%%%%%%%%%%%%%%%%%%%%%%%%%%%%%%%%%%%
%
%
%
\begin{The}\label{hqknegligibility}
Let $\Sigma\subset\Omega$ be a $C^1$ submanifold of dimension $p$.
Then the set of characteristic points $C(\Sigma)$ is $\cS^{p+1}$-negligible.
\end{The}
\begin{Rem}{\rm In order to apply the negligibility result of
\cite{Mag5} one has to check that the notion of characteristic
point in arbitrary stratified groups coincides with our
definition stated in the Heisenberg group.
According to \cite{Mag5} a point $x\in\Sigma$ is characteristic if
\begin{eqnarray}\label{charactpap}
\dim\lls H_x\H^{n}\rls-\dim\lls T_x\Sigma\cap H_x\H^{n}\rls\leq k-1.	
\end{eqnarray}
If $x$ is characteristic according to Definition~\ref{defchar}, then
$$
\dim\lls T_x\Sigma\cap H_x\H^{n}\rls=p=2n+1-k
$$
and (\ref{charactpap}) holds. Conversely, if (\ref{charactpap}) holds,
then 
$$
p=\dim(T_x\Sigma)=2n-k+1\leq \dim\lls T_x\Sigma\cap H_x\H^{n}\rls,
$$
hence $T_x\Sigma\subset H_x\H^{n}$.
}
\end{Rem}

%
%
%
%%%%%%%%%%%%%%%%%%%%%%%%%%%%%%%%%%%%%%%%%%%%%%%%%%%%%%%%%%%%%%
%   INTEGRAL FORMULA FOR THE SPHERICAL HAUSDORFF MEASURE
%%%%%%%%%%%%%%%%%%%%%%%%%%%%%%%%%%%%%%%%%%%%%%%%%%%%%%%%%%%%%%
%
%
%
\begin{Cor}\label{corinv}
Let $\Sigma\subset\Omega$ be a $C^1$ submanifold of dimension $p$.
Then we have
\begin{eqnarray}\label{thef2}
\int_\Sigma \theta\lls\tau_{\Sigma,\cV}(x)\rls\,d\cS^{p+1}(x)=\int_\Sigma |\tau_{\Sigma,\cV}(x)|\,d\vol_p(x)	
\end{eqnarray}
\end{Cor}
{\sc Proof}.
We apply Theorem~2.10.17(2) and Theorem~2.10.18(1) of \cite{Fed},
hence from limit (\ref{keyblow}) and Theorem~\ref{hqknegligibility} the proof follows by a standard argument. $\Box$
\begin{Rem}{\rm
Proof of Theorem~\ref{constmetintsp} immediately follows from (\ref{thef2}).
}\end{Rem}
Next, we present a class of homogeneous distances in the Heisenberg
group which possess constant metric factor.
The standard system of graded coordinates $F:\R^{2n+1}\lra\H^{n}$ induced by $(X_1,\ldots,X_{2n},Z)$ will be understood in the sequel.
To simplify notation we will write
$x=F(\xtil,x_{2n+1})\in\H^{n}$, with $\xtil=(x_1,\ldots,x_{2n})\in\R^{2n}$.
\begin{Pro}\label{homconstmet}
Let $F:\R^{2n+1}\lra\H^n$ define standard coordintates and let $\rho$ be a homogeneous distance of $\H^{n}$ such that
$\rho\lls 0,F(\cdot)\rls:\R^{2n+1}\lra\R$ only depends on
$(|\xtil|,x_{2n+1})$. Then $\theta_p^\rho(\tau)=\theta_p^\rho(\tilde{\tau})$
whenever $\tau,\tilde{\tau}$ are vertical simple $p$-vectors.
\end{Pro}
{\sc Proof}.
Let $\tau=U_1\wedge\cdots\wedge U_{p-1}\wedge Z$ and
$\tilde{\tau}=W_1\wedge\cdots\wedge W_{p-1}\wedge Z$ be vertical
simple $p$-vectors, where it is not restrictive assuming that both $(U_1,\cdots,U_{p-1},Z)$ and $(W_1,\ldots,W_{p-1},Z)$
are orthonormal systems of $\ch^{2n+1}$.
Then we easily find an isometry $J:\ch^{2n+1}\lra\ch^{2n+1}$
such that $J\lls\cL(\tau)\rls=\cL(\tilde{\tau})$ and $J(Z)=Z$.
Recall that our graded coordinates are defined by
$F=\exp\circ \cI$, where $\cI:\R^{2n+1}\lra\ch^{2n+1}$
is an isometry such that
$$
\cI(x_1,\ldots,x_{2n+1})=x_{2n+1}\,Z+\sum_{j=1}^{2n}x_j\,X_j
$$
for every $(x_1,\ldots,x_{2n+1})\in\R^{2n+1}$.
Thus, defining $\Btil_1=F^{-1}(B_1)\subset\R^{2n+1}$, we have
\begin{eqnarray}\label{rotcoord}
F^{-1}\lls \exp \cL(\tilde{\tau})\cap B_1\rls=
\cI^{-1}\circ J\lls \cL(\tau)\rls\cap\Btil_1=
\ph\lLs\cI^{-1}\lls \cL(\tau)\rls\rLs\cap\Btil_1\,,
\end{eqnarray}
where $\ph=\cI^{-1}\circ J\circ\cI:\R^{2n+1}\lra\R^{2n+1}$ is an Euclidean isometry such that $\ph(e_{2n+1})=e_{2n+1}$ and $e_{2n+1}$ is the
$(2n$$+$$1)$-th vector of the canonical basis of $\R^{2n+1}$.
Then $|\xtil|=|\ytil|$ whenever $\ph(\xtil,t)=(\ytil,t)$.
As a result, the fact that $\rho\lls 0,F(\xtil,t)\rls$ only depends on $(|\tilde{x}|,t)$
easily implies that $\ph(\Btil_1)=\Btil_1$. Thus, due to (\ref{rotcoord}),
it follows that
\begin{eqnarray}
\theta_p^\rho(\tau)=\mcH^p_{|\cdot|}\lLs\cI^{-1}\lls \cL(\tau)\rls
\cap\Btil_1\rLs=\mcH^p_{|\cdot|}\lLs\ph\lLs\cI^{-1}\lls \cL(\tau)\rls
\rLs\cap\Btil_1\rLs=\theta_p^\rho(\tilde{\tau}).
\end{eqnarray}
This ends the proof. $\Box$
%
%
%
%%%%%%%%%%%%%%%%%%%%%%%%%%%%%%%%%%%%%%%%%%%%%%%%%%%%%%%%%%%%%%
%EXAMPLES OF HOMOGENEOUS DISTANCES WITH CONSTANT METRIC FACTOR
%%%%%%%%%%%%%%%%%%%%%%%%%%%%%%%%%%%%%%%%%%%%%%%%%%%%%%%%%%%%%%
%
%
%
\begin{Exa}\label{examphomd}{\rm
An example of homogeneous distance satisfying hypotheses
of Proposition~\ref{homconstmet} is the gauge distance, also called
Kor\'anyi distance, \cite{Kor}.
The gauge distance from $x$ to the origin is given by
$$
d(x,0)=\lls|\xtil|^4+16\,x_{2n+1}^2\rls^{1/4},
$$
where $x=(\xtil,x_{2n+1})$. Then we define $d(x,y)=d(0,x^{-1}y)$, for any $x,y\in\H^{n}$.
Ano\-ther example of homogeneous distance with this property
is the ``maximum distance'', defined by
$$
d_\infty(x,0)=\max\left\{|\xtil|,|x_{2n+1}|^{1/2}\right\}.
$$
Due to Proposition~\ref{homconstmet}, both of these distances
have constant metric factor.
}\end{Exa}
\begin{Rem}{\rm
The metric factor depends on the Riemannian metric
$g$ we have fixed. Furthermore, if we divide it
by the volume of the unit ball $B_1$ (or any other fixed subset
of positive measure), then we obtain a number only depending
on the restriction of the Riemannian metric to $H\H^n$.
}\end{Rem} 
\begin{Lem}\label{labinv}
Let $\gtil$ be a left invariant Riemannian metric such that
$\gtil_{|H\H^n}=g_{|H\H^n}$. Then 
$
\theta_p^\rho(\tau)/\vol_p(B_1)=\tilde{\theta}_p^\rho(\tau)/
\widetilde{\vol}_p(B_1)
$
for any simple vertical $p$-vector, where $\widetilde{\vol}_p$
and $\tilde{\theta}_p^\rho(\tau)$ are defined with respect to the
metric $\gtil$.
\end{Lem}
{\sc Proof.}
By hypothesis, we can choose an orthonormal frame
$(X_1,X_2,\ldots,X_{2n},W)$ with respect to $\gtil$,
where $W=\lb\,Z+\sum_{j=1}^{2n}a_jX_j$ and $\lb\neq0$.
Let $F,\Ftil:\R^{2n+1}\lra\H^n$ represent system of
coordinates with respect to the standard basis
$(X_1,X_2,\ldots,X_{2n},Z)$ and $(X_1,X_2,\ldots,X_{2n},W)$,
respectively. We have $\Ftil=F\circ T$, where $T:\R^{2n+1}\lra\R^{2n+1}$
is given by the matrix
\begin{eqnarray}\label{matrixchalb}
A=\left[\begin{array}{cccccc}
1      & 0      & 0      & \cdots &  0       & a_1 \\
0      & 1      & 0      & \ddots &  0       & a_2 \\
\vdots & \ddots & 1      & \ddots &  \vdots  & a_3 \\
\vdots & \ddots & 0      & \ddots &  0       &\vdots \\
0      & \cdots & \cdots &  0 &  1       &a_{2n} \\  
0      & 0      & \cdots & 0      &  0       & \lb
 \end{array}\right]	
\end{eqnarray}
By Proposition~\ref{egriem}, it follows that
\begin{eqnarray}\label{lbvolprop}
\widetilde{\vol}_p=\Ftil_\sharp\cL^{2n+1}=|\det A|^{-1}F_\sharp\cL^{2n+1}
=|\lb|^{-1}\,\vol_p.
\end{eqnarray}
Let $\tau=U_1\wedge\cdots\wedge U_{p-1}\wedge Z$, where
it is not restrictive to assume that the horizontal
vectors $U_1,\ldots,U_{p-1}$ are orthonormal with respect to
both $g$ and $\gtil$. 
We denote by $\cL$ the subspace $\span\{U_1,\ldots,U_{p-1},Z\}$
of $\ch^n$. 
Recall that
\begin{eqnarray}
\tilde{\theta}_p^\rho(\tau)=\mcH^p_{|\cdot|}\lls\Ftil^{-1}(\exp(\cL)\cap B_1)\rls=
\mcH^p_{|\cdot|}\left( T^{-1}\lls F^{-1}(\exp(\cL)\cap B_1)\rls\right).
\end{eqnarray}
Now we wish to determine a basis of the subspace
$F^{-1}(\exp(\cL)\subset\R^{2n+1}$.
The relations $U_i=\sum_{j=1}^{2n}c_i^j\,X_j$ give rise to
$p-1$ orthonormal vectors $c_1,\ldots,c_{p-1}\in\R^{2n}\times\{0\}$
with respect to the Euclidean scalar product such that
\begin{eqnarray}
\span\{c_1,\ldots,c_{p-1},e_{2n+1}\}=F^{-1}
\left(\exp\lls\cL\rls\right)\subset\R^{2n+1}.
\end{eqnarray}
In order to compute the Euclidean jacobian of $T^{-1}$
restricted to the $p$-dimensional subspace $\span\{c_1,\ldots,c_{p-1},e_{2n+1}\}\subset\R^{2n+1}$
we write the matrix
\begin{eqnarray}\label{matrixchalb-1}
A^{-1}=\left[\begin{array}{cccccc}
1      & 0      & 0      & \cdots &  0       & -\lb^{-1}a_1 \\
0      & 1      & 0      & \ddots &  0       & -\lb^{-1}a_2 \\
\vdots & \ddots & 1      & \ddots &  \vdots  & -\lb^{-1}a_3 \\
\vdots & \ddots & 0      & \ddots &  0       &\vdots \\
0      & \cdots & \cdots &  0 &  1           &-\lb^{-1}a_{2n} \\  
0      & 0      & \ddots & \ddots &  0       & \lb^{-1}
\end{array}\right]	
\end{eqnarray}
noting that $T^{-1}c_j=c_j$ for every $j=1,\ldots,p-1$ and
$T^{-1}(e_{2n+1})=\lb^{-1}\,e_{2n+1}$.
As a result, the jacobian of
$\lls T^{-1}\rls_|:\span\{c_1,\ldots,c_{p-1},e_{2n+1}\}\lra\R^{2n+1}$
is $|\lb|^{-1}$, hence
\begin{eqnarray*}
&&\tilde{\theta}_p^\rho(\tau)=\mcH^p_{|\cdot|}\left( T^{-1}\lls F^{-1}(\exp(\cL)\cap B_1)\rls\right)\\
&&=|\lb|^{-1}\;\mcH^p_{|\cdot|}\lls(F^{-1}(\exp(\cL)\cap B_1)\rls
=|\lb|^{-1}\,\theta_p^\rho(\tau).
\end{eqnarray*}
Joining (\ref{lbvolprop}) and the previous equalities, our claim follows. $\Box$.
\begin{The}\label{invarresc}
Let $\Sigma$ be a $p$-dimensional submanifold of $\Omega$ and
let $\gtil$ a left invariant metric such that
$\gtil_{|H\H^n}=g_{|H\H^n}$. Then
\begin{eqnarray}
\frac{1}{\widetilde{\vol}_p(B_1)}\int_\Sigma |\tilde{\tau}_{\Sigma,\cV}(x)|\,d\widetilde{\vol}_p(x)
=\frac{1}{\vol_p(B_1)}\int_\Sigma |\tau_{\Sigma,\cV}(x)|\,d\vol_p(x)	
\end{eqnarray}
\end{The}
The previous theorem is an immediate conseqauence of Corollary~\ref{corinv} and Lemma~\ref{labinv}.
Next, we apply (\ref{thef2c}) to compute the spherical Hausdorff measure
of some submanifolds.
We will use the following proposition.
\begin{Pro}\label{concreteex}
Let $\phi:U\lra\R^{2n+1}$ be a $C^1$ embedding, where $U\subset\R^p$
is a bounded open set. Let $F:\R^{2n+1}\lra\H^n$ define standard coordinates and set $\Phi=F\circ\phi:U\lra\H^n$, where $\Sigma=\Phi(U)$.
Let $\rho$ be a homogeneous distance with constant
metric factor $\alpha>0$. Then we have
\begin{eqnarray}\label{auxforsp}
\cS^{p+1}_{\H^n}(\Sigma)=\int_U\left|\pi_\cV\lls\Phi_{u_1}(u)\wedge\Phi_{u_2}(u)\wedge\cdots\wedge\Phi_{u_p}(u)\rls\right|\,du\,,
\end{eqnarray}
for every measurable set $A\subset\H^n$, where the norm
$|\cdot|$ is induced by the scalar product (\ref{pscalarproduct})
on $p$-vectors.
\end{Pro}
{\sc Proof.}
By definition of Riemannian volume, formula (\ref{thef2c})
can be written with respect to $\phi$ as
\begin{eqnarray}\label{spexp}
\cS^{p+1}_{\H^n}(\Sigma)=\int_U|\tau_{\Sigma,\cV}(\phi(u))|\;
\sqrt{\det\left[g(\phi(u))\lls\Phi_{u_i}(u),\Phi_{u_j}(u)\rls\right]}
\;du\,,
\end{eqnarray}
where we have
\begin{eqnarray}\label{jacobwedge}
\left|\Phi_{u_1}(u)\wedge\Phi_{u_2}(u)\wedge\cdots\wedge\Phi_{u_p}(u)\right|=\sqrt{\det\left[g(\Phi(u))\lls\Phi_{u_i}(u),\Phi_{u_j}(u)\rls\right]}.
\end{eqnarray}
Therefore, taking into account the formula
\begin{eqnarray}\label{tauvert}
\tau_\Sigma\lls\Phi(u)\rls=\frac{\Phi_{u_1}(u)\wedge\Phi_{u_2}(u)\wedge\cdots\wedge\Phi_{u_p}(u)}{|\Phi_{u_1}(u)\wedge\Phi_{u_2}(u)\wedge\cdots\wedge\Phi_{u_p}(u)|}\,,
\end{eqnarray}
the definition of vertical tangent $p$-vector 
$\pi_\cV(\tau_\Sigma)=\tau_{\Sigma,\cV}$ and joining
(\ref{spexp}), (\ref{jacobwedge}) and (\ref{tauvert}),
formula (\ref{auxforsp}) follows. $\Box$

\begin{Exa}{\rm
Let $\phi:\R^3\lra\R^5$, defined by
$\phi(u)=(u_1,u_2,u_3,0,\frac{u_1^2+u_2^2+u_3^2}{2})$.
The mapping $\phi$ parametrizes a 3-dimensional
paraboloid $\Sigma=\Phi(U)$ of $\R^5$, where $U$ is an open bounded set of 
$\R^3$, $\Phi=F\circ\phi$ and $F:\R^3\lra\H^2$ represents
standard coordinates, according to Definition~\ref{gradef}.
Using expressions (\ref{hvf}), we have 
\begin{eqnarray*}
&&\phi_{u_1}(u)=\tilde{X}_1(\phi(u))+\lls\phi_3(u)+u_1\rls \tilde{T}(\phi(u))\,,\\
&&\phi_{u_2}(u)=\tilde{X}_2(\phi(u))+\lls\phi_4(u)+u_2\rls \tilde{T}(\phi(u))\,,   \\
&&\phi_{u_3}(u)=\tilde{X}_3(\phi(u))+\lls u_3-\phi_1(u)\rls \tilde{T}(\phi(u))\,.
\end{eqnarray*}
Observing that for every $j=1,\ldots,2n$, we have
$$
dF(\phi(u))\tilde{X_j}(\phi(u))=X_j(\Phi(u))\in H_{\Phi(u)}\H^n
\quad\mbox{and}\quad
$$
$$
dF(\phi(u))\tilde{Z}(\phi(u))=Z(\Phi(u))\in T_{\Phi(u)}\H^n\,,
$$
hence we obtain
\begin{eqnarray*}
&&\Phi_{u_1}(u)=X_1(\Phi(u))+\lls u_3+u_1\rls T(\Phi(u))\,,\quad
\Phi_{u_2}(u)=X_2(\Phi(u))+u_2T(\Phi(u))\,,   \\
&&\Phi_{u_3}(u)=X_3(\Phi(u))+\lls u_3-u_1\rls T(\Phi(u))\,.
\end{eqnarray*}
Thus, we can compute
\begin{eqnarray*}
&&\pi_\cV\lls\Phi_{u_1}\wedge\Phi_{u_2}\wedge\Phi_{u_3}\rls\\
&&=(u_3-u_1)\,X_1\wedge X_2\wedge T-u_2X_1\wedge X_3\wedge T
+(u_3+u_1)\,X_2\wedge X_3\wedge T\,,
\end{eqnarray*}
hence formula (\ref{auxforsp}) yields
\begin{eqnarray*}
\cS^4_{\H^2}(\Sigma)=\int_U	\sqrt{u_2^2+2(u_3^2+u_1^2)}\;du\,.
\end{eqnarray*}
}\end{Exa}
\begin{Exa}\label{compsppl}{\rm
Let $\phi:\R^2\lra\R^3$, $\phi(u_1,u_2)=(a_1u_1,a_2u_2,bu_1+cu_2)$, define a hyperplane in $\R^3$, where $a_1,a_2,b,c\in\R$ and
$\lls J\phi\rls^2=a_1^2a_2^2+a_1^2c^2+a_2^2b^2>0$.
Embedding the hyperplane in $\H^1$ through standard coordinates
$F:\R^3\lra\H^1$, we obtain
\begin{eqnarray*}
&&\Phi_{u_1}(u)=a_1X_1(\Phi(u))+(a_1a_2u_2+b)T(\Phi(u))\,,\\
&&\Phi_{u_2}(u)=a_2X_2(\Phi(u))+(c-a_1a_2u_1)T(\Phi(u))\,,
\end{eqnarray*}
where $\Phi=F\circ\phi$. Then we get
\begin{eqnarray*}
\pi_{\cV}\lls\Phi_{u_1}(u)\wedge\Phi_{u_2}(u)\rls=	
a_1(c-a_1a_2u_1) X_1\wedge T-(a_1a_2u_2+b)a_2X_2\wedge T
\end{eqnarray*}
and formula (\ref{auxforsp}) yields
\begin{eqnarray}\label{sphypplane}
\cS^3_{\H^1}(\Pi)=\int_U	\sqrt{a_1^2(c-a_1a_2u_1)^2+a_2^2(a_1a_2u_2+b)^2}\,du,
\end{eqnarray}
where $\Pi=\Phi(U)$ and $U$ is an open bounded set of $\R^2$.
}\end{Exa}
\begin{Exa}\label{compsppar}{\rm
Let $\phi:\R^2\lra\R^3$, $\phi(u_1,u_2)=(u_1,u_2,\frac{u_1^2+u_2^2}{2})$, define a paraboloid in $\R^3$. By standard coordinates $F:\R^3\lra\H^1$
and arguing as in the previous examples, we have
$$
\Phi_{u_1}(u)=X_1(\Phi(u))+(u_2+u_1)T(\Phi(u))\quad
\mbox{and}\quad
$$
$$
\Phi_{u_2}(u)=X_2(\Phi(u))+(u_2-u_1)T(\Phi(u))\,,
$$
where $\Phi=F\circ\phi$. It follows that
\begin{eqnarray*}
\pi_{\cV}\lls\Phi_{u_1}(u)\wedge\Phi_{u_2}(u)\rls=	
(u_2-u_1)X_1\wedge T-(u_2+u_1)X_2\wedge T
\end{eqnarray*}
and formula (\ref{auxforsp}) yields
\begin{eqnarray}\label{spparab}
\cS^3_{\H^1}(\cP)=\int_U	\sqrt{2u_1^2+2u_2^2}\;du
\end{eqnarray}
where $\cP=\Phi(U)$ and $U$ is an open bounded set of $\R^2$.
}\end{Exa}
\begin{Rem}{\rm
It is curious to notice that the density of $\cS^3_{\H^1}$
restricted to the paraboloid $\cP$, computed in (\ref{spparab}), is proportional to the density of $\cS^3_{\H^1}$ restricted to the horizontal
projection of $\cP$ onto the plane
$F\lls\{(x_1,x_2,x_3)\mid x_3=0\}\rls\subset\H^1$, whose density is given by (\ref{sphypplane}) in the case $a_1=a_2=1$ and $b=c=0$.
}\end{Rem}
\begin{Exa}{\rm
From computations of Example~\ref{compsppl}, one can get
the 2-dimensional spherical Hausdorff measure of the line
$\Phi(t)=F(at,0,bt)$ defined on an interval $[\alpha,\beta]$.
We have
\begin{eqnarray}
\Phi'(t)=a X_1(\Phi(t))+bT(\Phi(t))\quad\mbox{and}\quad
\pi_\cV(\Phi'(t))=bT(\Phi(t))\,,
\end{eqnarray}
then defining the submanifold $\cL=\Phi([\alpha,\beta])$, 
the formula $\cS^2_{\H^1}(\cL)=|b|\,(\beta-\alpha)$ holds.
}\end{Exa}
Another consequence of (\ref{thef2c}) is the lower semicontinuity of
the spherical Hausdorff measure with respect to
weak convergence of regular currents.
To see this, it suffices to establish the following formula
\begin{eqnarray}\label{supintegralfor}
\cS^{p+1}_{\H^n}(\Sigma)=\sup_
{\omega\in\cF_c^p(\Omega)}
\int_\Sigma\lan\tau_{\Sigma,\cV},\omega\ran\,d\vol_p\,,
\end{eqnarray}
where $\cF_c^p(\Omega)$ is the space of smooth $p$-forms with compact support in $\Omega$ with $|\omega|\leq1$. 
The norm of $\omega$ is defined making
the standard frame of $p$-forms $(dx_1,dx_2,\ldots,dx_{2n},\tilde{\theta})$
orthonormal and extending this scalar product to $p$-forms
exactly as we have seen in formula (\ref{pscalarproduct}).
The 1-form $\tilde{\theta}$ is the so called {\em contact form}
\begin{eqnarray}
\tilde{\theta}=dx_{2n+1}+\sum_{j=1}^n	x_{j+n}\;dx_j-x_j\;dx_{j+n}\,
\end{eqnarray}
written in standard coordinates. Note that $(dx_1,dx_2,\ldots,dx_{2n},\tilde{\theta)}$ is the dual basis
of $(\tilde{X}_1,\ldots,\tilde{X}_{2n},\tilde{Z})$.
Formula (\ref{supintegralfor}) follows from (\ref{thef2c})
observing that
\begin{eqnarray}
\int_\Sigma|\tau_{\Sigma,\cV}|\,d\vol_p=\sup_{\omega\in\cF_c^p(\Omega)}
\int_\Sigma\lan\tau_{\Sigma,\cV},\omega\ran\,d\vol_p\,,	
\end{eqnarray}
as one can check by standard arguments. As a consequence of these
observations, we can establish the following proposition.
\begin{Pro}
Let $(\Sigma_m)$ be a sequence of $C^1$ submanifolds of $\Omega$
which weakly converges in the sense of currents to
the $C^1$ submanifold $\Sigma$. Then 
\begin{eqnarray}\label{semic}
\liminf_{m\lra\infty}\cS^{p+1}_{\H^n}(\Sigma_m)\geq\cS^{p+1}_{\H^n}(\Sigma).	
\end{eqnarray}
\end{Pro}
{\sc Proof}.
By hypothesis
\begin{eqnarray}
\int_{\Sigma_m}\lan\tau_{\Sigma_m,\cV},\omega\ran\,d\vol_p
\lra\int_\Sigma\lan\tau_{\Sigma,\cV},\omega\ran\,d\vol_p\,,
\end{eqnarray}
for every $\omega\in\cF_c^p(\Omega)$. Then (\ref{supintegralfor}) ends the proof. $\Box$
\begin{Rem}{\rm
It is clear the importance of (\ref{semic}) in studying versions of the Plateau problem with respect to the geometry of Heisenberg groups.
}\end{Rem}
Recall that the horizontal normal is the orthogonal projection of the normal
to $\nu(x)$ to $T_x\Sigma$ onto the horizontal subspace $H_x\H^{n}$.
In the next proposition we show that in codimension one
an explicit relationship can be established between vertical tangent $2n$-vector and horizontal normal $\nu_H$.
\begin{Pro}\label{cod1h}
Let $\Sigma$ be a $2n$-dimensional submanifold of class $C^1$ and let $\nu_H(x)$ a horizontal normal at $x\in\Sigma$. Then we have
$$\nu_H^j=(-1)^j\;\tau_{\Sigma,\cV}^j$$
where $\nu_H=\sum_{j=1}^{2n}\nu_H^j\,X_j$ and
$
\tau_{\Sigma,\cV}=\sum_{j=1}^{2n}\tau_{\Sigma,\cV}^j\;
X_1\wedge\cdots X_{j-1}\wedge X_{j+1}\wedge\cdots X_{2n}\wedge Z.
$
In particular, the equality
$
|\tau_{\Sigma,\cV}|=|\nu_H|
$
holds.
\end{Pro}
{\sc Proof}
Let $(t_1,t_2,\ldots,t_{2n})$ be an orthonormal basis of $T_x\Sigma$, where
$x$ is a transverse point. Then
$$
t_j=\sum_{i=1}^{2n}\,c_j^i\; X_i(x)+c_j^{2n+1} Z(x)
$$
where $C=(c_j^i)$ is a $(2n+1)\times 2n$ matrix, whose columns
are orthonormal vectors of $\R^{2n+1}$.
Then we have
\begin{eqnarray*}
\tau_\Sigma(x)=t_1\wedge t_2\wedge\cdots\wedge t_{2n}=\sum_{j=1}^{2n+1}
\det(\mbox{{\rm \^C}}^{\;j}) \;
X_{1}\wedge X_{2}\wedge\cdots\wedge X_{j-1}\wedge X_{j+1}\wedge\cdots\wedge Z\,,
\end{eqnarray*}
where $\mbox{{\rm \^C}}^{\;j}$ is the $2n\times2n$ matrix obtained by removing
the $j$-th row from $C$. The vertical projection yields
\begin{eqnarray}\label{vertangmat}
\tau_{\Sigma,\cV}(x)=\pi_\cV\lls\tau_\Sigma(x)\rls=\sum_{j=1}^{2n}
\det(\mbox{{\rm \^C}}^{\;j}) \;
X_{1}\wedge X_{2}\wedge\cdots\wedge X_{j-1}\wedge X_{j+1}\wedge\cdots\wedge Z
\end{eqnarray}
and by elementary linear algebra one can deduce that
\begin{eqnarray}
\sum_{j=1}^{2n+1} (-1)^{j}\det(\mbox{{\rm \^C}}^{\;j})\;c_k^j=
\det\left[\begin{array}{cc}
\!C\! & \!c_k\! \end{array}\right]=0
\end{eqnarray}
for every $k=1,\ldots,2n$. Then the vector
$$
\nu=\sum_{j=1}^{2n} (-1)^{j}\det(\mbox{{\rm \^C}}^{\;j})\;X_j
+(-1)^{2n+1}\det(\mbox{{\rm \^C}}^{\;2n+1})\;Z
$$
yields a unit normal to $\Sigma$ at $x$.
Its horizontal projection is
\begin{eqnarray}\label{hornormat}
\nu_H=\sum_{j=1}^{2n} (-1)^{j}\det(\mbox{{\rm \^C}}^{\;j})\;X_j.
\end{eqnarray}
Formulae (\ref{vertangmat}) and (\ref{hornormat}) yield the thesis. $\Box$
%
%
%
%
%%%%%%%%%%%%%%%%%%%%%%%%%%%%%%%%%%%%%%%%%%%%%%%%%%%%%%%%%%%%%
%
%             SECTION ON COAREA FORMULA 
%
%%%%%%%%%%%%%%%%%%%%%%%%%%%%%%%%%%%%%%%%%%%%%%%%%%%%%%%%%%%%%
%
%
%
%
\section{Coarea formula}\label{coareasect}
This section is devoted to the proof of Theorem~\ref{coacodk}.
Next, we recall the Riemannian coarea formula, see Section~13.4 of \cite{Bur}.
\begin{The}
Let $f:\H^{n}\lra\R^k$ be a Riemannian Lipschitz function,
with $1\leq k<2n+1$.
Then for any summable map $u:\H^{n}\lra\R$, the following formula holds
\begin{equation}\label{riem}
\int_{\H^{n}}\,u(x)\;J_gf(x)\,d\vol_{2n+1}(x)
=\int_{\R^k}\bigg(\int_{f^{-1}(t)}u(y)\,d\vol_p(y)\bigg)dt\,,
\end{equation}
where $p=2n+1-k$
\end{The}
In the previous theorem the Heisenberg group $\H^{n}$
is equipped with its left invariant Riemannian metric $g$.
The terminology ``Riemannian Lipschitz map" means that the map
is Lipschitz with respect to the Riemannian distance.
%
%
%
%%%%%%%%%%%%%%%%%%%%%%%%%%%%%%%%%%%%%%%%%%%%%%%%%%%%%%%%%%%%%%
%            PROOF OF COAREA FORMULA
%%%%%%%%%%%%%%%%%%%%%%%%%%%%%%%%%%%%%%%%%%%%%%%%%%%%%%%%%%%%%%
%
%
%
\vskip.25truecm
\noindent
{\em Proof of Theorem~\ref{coacodk}}.
We first prove (\ref{cht}) in the case $f$ is defined on all of $\H^{n}$
and is of class $C^1$. Let $\Omega$ be an open subset of $\H^{n}$.
In view of Riemannian coarea formula (\ref{riem}), we have
\begin{eqnarray}\label{eucl2}
\int_\Omega u(x) J_gf(x)\,dx=\int_{\R^k}\left(\int_{f^{-1}(t)\cap\Omega}
u(y)\,d\vol_p(y)\right)dt,
\end{eqnarray}
where $u:\Omega\lra[0,+\infty]$ is a measurable function.
Note that in the left hand side of (\ref{eucl2}) we have used
the Lebesgue measure in that, by Proposition~\ref{egriem},
it coincides with the volume measure
expressed in terms of standard coordinates, namely
$F_\sharp(\cL^{2n+1})=\vol_{2n+1}$.
Now we define
$$
u(x)=J_{H}f(x){\bf 
1}_{\{Jf\neq0\}\cap\Omega}(x)/Jf(x)
$$
and use (\ref{eucl2}), obtaining
\begin{eqnarray}\label{eucr2}
\int_{\Omega} J_{H}f(x)\,dx=\int_{\R^k}\left(\int_{f^{-1}(t)\cap\Omega}
\frac{J_{H}f(x){\bf 1}_{\{Jf\neq0\}}(x)}{Jf(x)}\;d\vol_p(y)\right)dt.
\end{eqnarray}
The validity of (\ref{eucl2}) also implies that for a.e. $t\in\R^{k}$ 
the set of points of $f^{-1}(t)$ where $J_gf$ vanishes is 
$\vol_p$-negligible, then the previous formula becomes
\begin{eqnarray}\label{eucrr2}
\int_{\Omega}J_{H}f(x)\,dx=\int_{\R^k}\left(\int_{f^{-1}(t)\cap\Omega}
\frac{J_{H}f(x)}{Jf(x)}\;d\vol_p(y)\right)dt.
\end{eqnarray}
By classical Sard's theorem and Theorem~\ref{hqknegligibility}
for a.e. $t\in\R^{k}$ the $C^1$ submanifold $f^{-1}(t)$
has $\cS_{\H^n}^{p+1}$-negligible characteristic points, hence
Proposition~\ref{charcharhj} implies that
$$
C_{t}=\{y\in f^{-1}(t)\cap\Omega\mid J_{H}f(y)=0\}
$$
is $\cS_{\H^n}^{p+1}$-negligible. As a result, from formulae (\ref{jhoverjg}) and (\ref{thef2c})
we have proved that for a.e. $t\in\R^{k}$ the equalities
$$
\int_{f^{-1}(t)\cap\Omega}
\frac{J_{H}f(x)}{J_gf(x)}\;d\vol_p(y)=\cS_{\H^{n}}^{p+1}(f^{-1}(t)
\cap\Omega\sm C_{t})=\cS^{p+1}_{\H^{n}}(f^{-1}(t)\cap\Omega)
$$
hold, therefore (\ref{eucrr2}) yields
\begin{eqnarray}\label{chh2}
\int_{\Omega}J_{H}f(x)\,dx
=\int_{\R^k}\cS^{p+1}_{\H^{n}}(f^{-1}(t)\cap\Omega)\;dt.
\end{eqnarray}
The arbitrary choice of $\Omega$ yields the validity of
(\ref{chh2}) also for arbitrary closed sets.
Then, approximation of measurable sets by closed ones,
Borel regularity of $\cS^{p+1}_{\H^{n}}$
and the coarea estimate 2.10.25 of \cite{Fed}
extend the validity of (\ref{chh2}) to the following one
\begin{eqnarray}\label{chh}
\int_{A}J_{H}f(x)\,dx=\int_{\R^k}\cS^{p+1}_{\H^n}(f^{-1}(t)\cap A)dt,
\end{eqnarray}
where $A$ is a measurable subset of $\H^{n}$.
Now we consider the general case, where $f:A\lra\R^k$ is a Lipschitz
map defined on a measurable bounded subset $A$ of $\H^3$.
Let $f_1:\H^{n}\lra\R^k$ be a Lipschitz extension of $f$, namely,
${f_1}_{|A}=f$ holds.
Due to the Whitney extension theorem (see for instance 3.1.15 of \cite{Fed})
for every arbitrarily fixed $\ep>0$ there exists
a $C^1$ function $f_2:\H^{n}\lra\R^k$ such that the open subset
$O=\{z\in\H^{n}\mid f_1(z)\neq f_2(z)\}$ has Lebesgue measure
less than or equal to $\ep$.
We wish to prove 
\begin{eqnarray}\label{cancon}
&&\left|\int_A J_Hf(x)\,dx-\int_{\R^k}\cS^{p+1}_{\H^n}(f^{-1}(t)\cap A)dt
\right|\leq \int_{A\cap O} J_Hf(x)\,dx \nonumber \\ 
&&+\int_{\R^k}\cS^{p+1}_{\H^n}(f^{-1}(t)\cap A\cap O)dt\,.
\end{eqnarray}
In fact, due to the validity of (\ref{chh}) for $C^1$ mappings, we have
$$
\int_{A\sm O} J_Hf_2(x)\,dx
=\int_{\R^k}\cS^{p+1}_{\H^{n}}(f_2^{-1}(t)\cap A\sm O)\,dt.
$$
Note here that the horizontal jacobian $J_Hf$ is well defined
on $A$, in that $df$ is well defined at density points of the domain, see for instance Definition~7 and Proposition~2.2 of \cite{Mag}.
The equality ${f_2}_{|A\sm O}=f_{|A\sm O}$
implies that $J_Hf_2=J_Hf$ a.e. on $A\sm O$, therefore
$$
\int_{A\sm O} J_Hf(x)\,dx
=\int_{\R^k}\cS^{p+1}_{\H^{n}}(f^{-1}(t)\cap A\sm O)\,dt
$$
holds and inequality (\ref{cancon}) is proved.
Now we observe that for a.e. $x\in A$, we have
$$
J_Hf(x)\leq \prod_{i=1}^k\bigg(\sum_{j=1}^{2n}\lLs X_jf^i(x)\rLs^2\bigg)^{1/2}
\leq\|df(x)_{|H_x\H^n}\|^k
$$
therefore the estimate
\begin{eqnarray}\label{efst}
J_Hf(x)\leq \Lip(f)^k
\end{eqnarray}
holds for a.e. $x\in A$. By virtue of the general coarea inequality
2.10.25 of \cite{Fed} there exists a dimensional constant $c_1>0$ 
such that
\begin{eqnarray}
\int_{\R^k}\cS^{p+1}_{\H^{n}}(f^{-1}(t)\cap A\cap O)dt\leq c_1\;
\Lip(f)^k\;\mcH^{2n+2}(O).
\end{eqnarray}
The fact that the $2n+2$-dimensional Hausdorff measure $\mcH^{2n+2}$
with respect to the homogeneous distance $\rho$ is proportional
to the Lebesgue measure, gives us a constant $c_2>0$ such that
\begin{eqnarray}\label{labelineq}
\int_{\R^k}\cS^{p+1}_{\H^{n}}(f^{-1}(t)\cap A\cap O)dt\leq c_2\;
\Lip(f)^k\;\cL^{2n+1}(O)\leq  c_2\;\Lip(f)^k\;\ep.
\end{eqnarray}
Thus, estimates (\ref{efst}) and (\ref{labelineq}) joined with inequality
(\ref{cancon}) yield 
\begin{eqnarray*}
\left|\int_A J_Hf(x)\,dx-\int_{\R^k}\cS^{p+1}_{\H^{n}}(f^{-1}(t)\cap A)dt
\right|\leq \,(1+c_2)\;\Lip(f)^k\;\ep.
\end{eqnarray*}
Letting $\ep\ra0^+$, we have proved that
\begin{eqnarray}
\int_A J_Hf(x)\,dx=\int_{\R^k}\cS^{p+1}_{\H^{n}}(f^{-1}(t)\cap A)dt.
\end{eqnarray}
Finally, utilizing increasing sequences of step functions
pointwise converging to $u$ and applying Beppo Levi
convergence theorem the proof of (\ref{cht}) is achieved
in the case $A$ is bounded. If $A$ is not bounded, then one
can take the limit of (\ref{cht}) where $A$ is replaced by $A_k$ and
$\{A_k\}$ is an increasing sequence of measurable bounded sets
whose union yields $A$. Then the Beppo Levi
convergence theorem concludes the proof. $\Box$
\begin{Rem}{\rm
Notice that once $f:A\lra\R^k$ in the previous theorem
is considered with respect to standard coordinates it is
easy to check that the locally Lipschitz property with respect to the Euclidean distance of $\R^{2n+1}$ is equivalent to the locally
Lipschitz property with respect to the Riemannian distance.
}\end{Rem}

\end{document}